\documentclass[12pt,letterpaper]{article}
\pdfoutput=1

\usepackage{graphics,epsfig,graphicx,color}
\graphicspath{{/EPSF/}{../Figures/}{Figures/}}

\usepackage{amssymb,latexsym}

\usepackage{setspace}

\usepackage{amsmath}

\usepackage{mathrsfs,amsfonts}

\usepackage{amsthm,amsxtra}

\newif\ifPDF
\ifx\pdfoutput\undefined
	\PDFfalse
\else
	\ifnum\pdfoutput > 0
		\PDFtrue
	\else
		\PDFfalse
	\fi
\fi

\ifPDF
	\usepackage{pdftricks}
	\begin{psinputs}
		\usepackage{pstricks}
		\usepackage{pstcol}
		\usepackage{pst-plot}
		\usepackage{pst-tree}
		\usepackage{pst-eps}
		\usepackage{multido}
		\usepackage{pst-node}
		\usepackage{pst-eps}
	\end{psinputs}
\else
	\usepackage{pstricks}
\fi

\ifPDF
	\usepackage[debug,pdftex,colorlinks=true, 
	linkcolor=blue, bookmarksopen=false,
	plainpages=false,pdfpagelabels]{hyperref}
\else
	\usepackage[dvips]{hyperref}
\fi

\usepackage{fancyhdr}

\usepackage{comment}

\newcommand{\dint}{\displaystyle\int}

\newcommand{\bPhi}{\boldsymbol \Phi}
\newcommand{\bPsi}{\boldsymbol \Psi}

\newcommand{\bSigma}{{\boldsymbol\Sigma}}
\newcommand{\bLambda}{{\boldsymbol\Lambda}}

 \newcommand{\bbN}{\mathbb N}
 
\newcommand{\bbR}{\mathbb R} \newcommand{\bbS}{\mathbb S}

 \newcommand{\bn}{\mathbf n}

 \newcommand{\bv}{\mathbf v} 
 \newcommand{\bx}{\mathbf x} 
\newcommand{\by}{\mathbf y} \newcommand{\bz}{\mathbf z} 
\newcommand{\bA}{\mathbf A} \newcommand{\bB}{\mathbf B}

\newcommand{\bG}{\mathbf G} \newcommand{\bH}{\mathbf H}
\newcommand{\bI}{\mathbf I} 
\newcommand{\bK}{\mathbf K}

\newcommand{\bS}{\mathbf S}

\newcommand{\cA}{\mathcal A} \newcommand{\cB}{\mathcal B}

\newcommand{\cG}{\mathcal G} 
 
\newcommand{\cK}{\mathcal K} 
\newcommand{\cM}{\mathcal M} 
\newcommand{\cO}{\mathcal O}

\newcommand{\fA}{\mathfrak A} \newcommand{\fB}{\mathfrak B}
\newcommand{\fE}{\mathfrak E} 
\newcommand{\fU}{\mathfrak U}
\newcommand{\fa}{\mathfrak a} 
\newcommand{\ff}{\mathfrak f}
\newcommand{\fs}{\mathfrak s}
\newcommand{\ft}{\mathfrak t}
\newcommand{\fu}{\mathfrak u}
\newcommand{\fx}{\mathfrak x}

\newcommand{\wU}{{\widetilde U}}
\newcommand{\rd}{{\rm d}}

\DeclareMathOperator*{\argmin}{arg\,min}

\newenvironment{keywords}
{\noindent{\bf Key words.}\small}{\par\vspace{1ex}}

\title{Inverse transport calculations in optical imaging with subspace optimization algorithms}

\author{
	Tian Ding\thanks{
		Department of Mathematics,
		University of Texas, 
		Austin, TX 78712; tding@math.utexas.edu
	}
	\and 
	Kui Ren\thanks{
		Department of Mathematics,
		University of Texas, 
		Austin, TX 78712; ren@math.utexas.edu
	}
}

\begin{document}

\maketitle



\begin{abstract}
	Inverse boundary value problems for the radiative transport equation play important roles in optics-based medical imaging techniques such as diffuse optical tomography (DOT) and fluorescence optical tomography(FOT). Despite the rapid progress in the mathematical theory and numerical computation of these inverse problems in recent years, developing robust and efficient reconstruction algorithms remains as a challenging task and an active research topic. We propose here a robust reconstruction method that is based on subspace minimization techniques. The method splits the unknown transport solution (or a functional of it) into low-frequency and high-frequency components, and uses singular value decomposition to analytically recover part of low-frequency information. Minimization is then applied to recover part of the high-frequency components of the unknowns. We present some numerical simulations with synthetic data to demonstrate the performance of the proposed algorithm.
\end{abstract}






\section{Introduction}
\label{SEC:Intro}

The mathematical and computational study of inverse coefficient problems to the radiative transport equations have attracted extensive attentions in recent years; see for instance the reviews~\cite{Bal-IP09,McCormick-TTSP84,McCormick-TTSP86,McCormick-NSE92,Ren-CiCP10} and references therein. The main objective of these inverse problems is to reconstruct physical parameters in the radiative transport equation from partial information on the solution to the equation. These inverse problems have important applications in many areas of science and engineering, such as ocean, atmospheric and interstellar optics~\cite{Belleni-Morante-TTSP03,BaRe-IP05,Chahine-JAS70,ElDuMcEm-JOSA88,Li-JQSRT01,McCormick-JOSA04,SoVeStCh-TTSP04,WaUe-ASS89}, radiation therapy planning~\cite{ArBaLuGr-MMMTT91,Boman-Thesis07,GiHoWaFaMo-PMB06,JiScPaJi-PMB12,Luo-RPC98,ShFeOlMa-SIAM99,TiHeToSiAlPyUl-PMB08,VaWaMcFaSaMo-PMB10}, diffuse optical tomography and quantitative photoacoustic tomography~\cite{AbReHi-IP05,ArBaLuGr-MMMTT91,Arridge-IP99,Arridge-IP06,BaTa-SIAM07,CaXuAl-IEEE03,Dorn-CAMQ01,GaZh-TTSP09,GoKi-IP09,GuReHi-AO07,KiFlYaKaHi-BOE10,KiHi-IP09,KlYa-SIAM07,Langmore-IP08,MaRe-CMS14,ReBaHi-SIAM06,ReBaHi-AO07,ScMa-IPI07,Tamasan-IP02,Tamasan-CM03,TaHaHa-IP13,TaVaAr-JQSRT08,Wang-AIHP99}, molecular imaging~\cite{BaTa-SIAM07,GaZh-OE10A,GaZh-OE10B,KiMo-IP06,KlNtHi-JCP05,StUh-APDE08} and many more~\cite{BeRo-JMAA02,Babovsky-SIAM91,Babovsky-IP95,Barichello-Pro02,BeRo-JMAA02,Case-PF73,Case-PF75,Case-PF77,ChSt-CPDE96,ChSt-IP96,DrBa-M2AS05,KaDa-TTSP79,KaMo-JMP78,KaMo-JMP78B,KiIs-JCP99,KlYa-SIAM07,Larsen-TTSP88,Panchenko-IP93,SoVeStCh-TTSP04,TiRoNeCa-Pro02,Wang-AMC00,WaUe-ASS89,YiWuSu-JCP96}.

We consider here the application of inverse transport problems in biomedical optical imaging techniques such as diffuse optical tomography (DOT)~\cite{AbReHi-IP05,ArBaLuGr-MMMTT91,Arridge-IP99,Arridge-IP06,BaTa-SIAM07,CaXuAl-IEEE03,Dorn-CAMQ01,GaZh-TTSP09,GoKi-IP09,GuReHi-AO07,KlYa-SIAM07,ReBaHi-SIAM06,ReBaHi-AO07,ScMa-IPI07} and fluorescence optical tomography (FOT)~\cite{BaTa-SIAM07,GaZh-OE10A,GaZh-OE10B,KlNtHi-JCP05} where radiative transport equations are often employed as the model for light propagation in biological tissues. To setup the problem, let us denote by $\Omega\subset\bbR^d$ ($d\ge 2$) the tissue of interest, with sufficiently regular surface $\partial\Omega$. We denote by $\bbS^{d-1}$ the unit sphere in $\bbR^d$, and $\bv \in\bbS^{d-1}$ the unit vector on the sphere. We denote by $X\equiv\Omega\times \bbS^{d-1}$ the phase space and define the boundary sets of the phase space, $\Gamma_{\pm}$, as
\[
	\Gamma_\pm= \{ (\bx,\bv) \in \partial\Omega\times \bbS^{d-1} \mbox{ s.t. } \pm\bv\cdot\bn(\bx)>0\}
\]
with $\bn(\bx)$ the unit outer normal vector at $\bx\in\partial\Omega$. The radiative transport equation for the phase-space density distribution of the photons in the tissue can be written as
\begin{equation}\label{EQ:ERT}
	\begin{array}{rcll}
  	\bv\cdot\nabla u(\bx,\bv) + \sigma_\fa(\bx)u(\bx,\bv) &=& \sigma_\fs(\bx)\cK u(\bx,\bv) &\mbox{in}\ \ X\\
    u(\bx,\bv) &=& f(\bx) & \mbox{on}\ \ \Gamma_{-},
	\end{array}
\end{equation}
where $u(\bx,\bv)$ is the density of photons at $\bx\in\Omega$ traveling in direction $\bv$, and $f$ is the light source. The positive functions $\sigma_\fa(\bx)$ and $\sigma_\fs(\bx)$ are the absorption coefficient and the scattering coefficients respectively. The total absorption coefficient is given by $\sigma(\bx)\equiv\sigma_\fa(\bx)+\sigma_\fs(\bx)$. The scattering operator $\cK$ is given by
\begin{equation}\label{EQ:Scat Operator}
	\cK u (\bx,\bv)=\dint_{\bbS^{d-1}}k(\bv, \bv')u(\bx,\bv')d\bv' - u(\bx,\bv)
\end{equation}
where the scattering kernel $k(\bv,\bv')$ describes the probability that photons traveling in direction $\bv'$ getting scattered into direction $\bv$. Note that to conserve the total mass, we have normalized the surface measure $d\bv$ on $\bbS^{d-1}$ and the scattering kernel $k(\bv,\bv')$ such that
\begin{equation}\label{EQ:Normalization}
	\int_{\bbS^{d-1}} d\bv=1, \qquad \mbox{and} \qquad 
	\int_{\bbS^{d-1}}k(\bv, \bv')d \bv'=1,\ \ \forall\bv\in \bbS^{d-1}.
\end{equation}
In biomedical optics, the scattering kernel is often taken as the Henyey-Greenstein phase function~\cite{HeGr-AJ41}:
\begin{equation}\label{EQ:HG}
k(\bv, \bv')\equiv k_g(\bv\cdot\bv')=\Pi\dfrac{1-g^2}{(1+g^2-2g\bv\cdot\bv')^{d/2}}, 
\end{equation}
which is a one-parameter function that depends only on the angle between the two directions $\bv$ and $\bv'$ for a given anisotropy factor $g\in[-1,1]$. The normalization constant $\Pi$ is determined by the normalization condition~\eqref{EQ:Normalization}.

The function $f(\bx)$ models the illumination source used in imaging experiments. In practical application of biomedical imaging, for instance in DOT and FOT, it is often technically difficult to construct angularly-resolved illumination sources. This is the main reason for us to employ an isotropic source function (referring to the fact that $f(\bx)$ does not depend on the angular variable $\bv$) in the the transport model~\eqref{EQ:ERT}. The measured data in biomedical optical imaging is usually a functional of the solution to the transport equation. Once again, due to the fact that it is difficult to measure angularly-resolved quantities, angularly-averaged quantities are usually measured. Here we consider applications (for instance DOT and FOT) where the measurement is the photon current on the surface of the tissue. The current is defined as
\begin{equation}\label{EQ:Data}
	j(\bx)\equiv \cM u(\bx)=\int_{\{\bv\in\bbS^{d-1}:\ \bv\cdot\bn(\bx)>0\}} \bv\cdot \bn(\bx) u(\bx,\bv)_{|\Gamma_+} d\bv, \quad \bx\in\partial\Omega.
\end{equation}
The objective of the biomedical imaging problems here is to reconstruct the optical absorption and scattering coefficients of biological tissues, $\sigma_\fa$ and $\sigma_\fs$ from data encoded in the albedo operator:
\begin{equation}\label{EQ:Albedo}
	\Lambda_{\sigma_\fa,\sigma_\fs}:\ \ f(\bx) \mapsto j(\bx)
\end{equation}
   
There are two major issues with diffuse optical imaging. The first issue is its low resolution due to the multiple scattering of light in tissues. Mathematically, this is manifested as the instability of the inverse transport problem~\cite{Bal-IP09,BaLaMo-IPI08}. By instability we mean that the noise in the data are significantly amplified in the inversion process, assuming that the problem admits a unique solution to start with. To stabilize the inverse problem, one can incorporate additional \emph{a priori} information into the computational inversion algorithms. Commonly-used \emph{a priori} information including, for instance, the smoothness or non-smoothness of the unknown~\cite{GaZh-OE10A,GaZh-OE10B,GuReMaHi-JBO10} and the shape of the regions of interests~\cite{Arridge-IP06,Dorn-CAMQ01}. The second issue with diffuse optical tomography is that there is no analytical inversion formulas for the image reconstruction problem, even in very academic geometrical configuration~\cite{ScMa-IPI07}. Computational reconstruction algorithms based on the radiative transport model are in general extremely slow. Fast reconstruction algorithms are actively sought by researchers in the field.

The instability of the inverse transport problems implies that when there is no available \emph{a priori} information, only low-frequency components of the unknowns can be reconstructed stably. One should thus not spend too much efforts trying to reconstruct high-frequency components of the unknowns. Based on this observation and the idea of subspace minimization~\cite{Chen-JOSA09,PaChZhYe-JOSA10,ZhCh-IP09}, we propose here a fast computational reconstruction method for the aforementioned inverse transport problems. Our method relies on the fact that we can explicitly factorize out some unstable components of the albedo operator $\Lambda_{\sigma_\fa,\sigma_\fs}$ defined in~\eqref{EQ:Albedo}. The unstable components of the albedo operator then impose a natural limit on the highest-frequency components of the unknown that can be reconstructed stably from the data. The factorization of $\Lambda_{\sigma_\fa,\sigma_\fs}$ is not unique in general. For our purpose, we follow the ideas in~\cite{Chen-JOSA09,PaChZhYe-JOSA10,ZhCh-IP09} to reformulate the transport problem into the form
\begin{eqnarray}
	\label{EQ:ERT SOM A} j &=& \cA \fu, \\
	\label{EQ:ERT SOM B} \fu &=& \cB \fu - \ff,
\end{eqnarray}
where $\cA$ is an operator that \emph{does not depend on the unknowns}, and $\cB$ is an operator that only depends on the unknowns. The intermediate quantity $\fu$ can be either the transport solution $u$, or a functional of it, or a functional  of both $u$ and the unknown coefficient, depending on the setup of the inverse problem. We will focus on the derivation of this formulation in the next two sections.

Once the forward problem is put into this form, we can split the inverse problem into two steps. In the first step, we ``invert'' the first equation to find $\fu$ while in the second step, we ``invert'' from $\fu$ the unknown coefficients. By construction, $\cA$ can not be inverted or not stably inverted. This means that there are many small or even zero singular values of $\cA$. Noise in the data can easily ruin the reconstruction of $\fu$ on these components of the small singular values. We thus split the eigenspace of $\cA$ into two subspaces, one spanned by eigenvectors corresponding to large singular values $\mu$ (say $\mu \ge \mu_c$) that we call signal subspace, and the other spanned by eigenvectors corresponding to small singular values ($\mu <\mu_c$) that we call noise subspace, following terminologies in~\cite{Chen-JOSA09,PaChZhYe-JOSA10,ZhCh-IP09}. The unknown coefficient functions and the intermediate quantity $\fu$ can be then be written as summations of two components that correspond to their projections into the two subspaces. This explicit splitting allows us to focus on the components of the unknowns corresponding to signal subspace while pay less attention to, and sometimes completely throw away, the components corresponding to noise subspace. The price we have to pay is that we have to be able to construct the singular value decomposition (SVD) of the operator $\cA$ in a computationally inexpensive way so that the whole algorithm is computationally feasible. Indeed, as we will see later, due to the fact that $\cA$ is independent of the unknowns, the SVD of $\cA$ can be pre-computed and do \emph{not} need to be updated in the nonlinear minimization process. This is the main reason why the reconstruction algorithm can be efficient.

The rest of the paper is organized as follows. In Section~\ref{SEC:Cont Form}, we reformulate, on continuous level, the radiative transport equation into the form of system~\eqref{EQ:ERT SOM A} and~\eqref{EQ:ERT SOM B} for both the reconstruction of the absorption coefficient and the reconstruction of the scattering coefficient. We then present in Section~\ref{SEC:Disc Form} the discretization of the continuous formulations. The details of our subspace-based minimization algorithm is then presented in Section~\ref{SEC:SOM} with some numerical experiments to demonstrate its performance in Section~\ref{SEC:Num}. Concluding remarks are then offered in Section~\ref{SEC:Concl}.

\section{Continuous Formulation}
\label{SEC:Cont Form}

It is well-known in inverse transport theory that with the types of data encoded in the albedo operator $\Lambda_{\sigma_\fa,\sigma_\fs}$ defined in~\eqref{EQ:Albedo}, only one of the two optical coefficients can be uniquely reconstructed~\cite{Bal-IP09,BaLaMo-IPI08} when no further information are available. We will thus work on the reconstruction of one coefficient assuming that the other is known.

\subsection{Recovering absorption coefficient}
\label{SUBSEC:Absorption}

Let us first consider the case of reconstructing the absorption coefficient $\sigma_\fa$ assuming that the scattering coefficient $\sigma_\fs$ is known. We first introduce the adjoint boundary Green function $G^\fa_b(\bx,\bv;\by)$ as the solution of the following adjoint free transport equation:
\begin{equation}\label{EQ:Green A B}
  \begin{array}{rcll}
  	-\bv\cdot\nabla G^\fa_b(\bx,\bv;\by) &=& 0, &\mbox{in}\ \ X\\
	G^\fa_b(\bx,\bv;\by) & =& \delta(\bx-\by), &\mbox{on}\ \ \Gamma_{+} .
  \end{array}
\end{equation}
Thus $G^\fa_b(\bx,\bv;\by)$ is the solution of the adjoint transport problem with an isotropic point source on the outgoing boundary $\Gamma_+$. Later on in this paper, $\by\in\partial\Omega$ is considered as the location of the detector used to measure the current data $j(\bx)$.

We now multiply~\eqref{EQ:ERT} by $G^\fa_b(\bx,\bv;\by)$, subtract it by the multiplication of ~\eqref{EQ:Green A B} and $u(\bx,\bv)$, and then integrate over phase space $X$, we have
\begin{equation}\label{EQ:Phase ID A B}
   \int_{X}\big(G^\fa_b \bv\cdot\nabla u + u\bv\cdot\nabla G^\fa_b\big) \rd\bx \rd\bv = \int_{X} G^\fa_b \big(\sigma_\fs \cK u  - \sigma_\fa  u\big) \rd\bx \rd\bv.
\end{equation}
The left hand side can be simplified using integration-by-part. We find, after splitting the boundary integral into an integral on $\Gamma_-$ and another integral on $\Gamma_+$, that
\begin{equation}\label{EQ:IBP A B}
	\int_{X}\big(G^\fa_b \bv\cdot\nabla u + u\bv\cdot\nabla G^\fa_b\big) \rd\bx \rd\bv = \cM u(\by) + \int_{\Gamma_-} \bv\cdot \bn G^\fa_b(\bx,\bv;\by) f(\bx) \rd S(\bx) \rd\bv,
\end{equation}
where $\rd S(\bx)$ is the surface measure on $\partial\Omega$.

We combine ~\eqref{EQ:Phase ID A B} and ~\eqref{EQ:IBP A B} to obtain
\begin{equation}\label{EQ:Data A B}
	j(\by)\equiv \cM u(\by) = \int_{X} G^\fa_b \big(\sigma_\fs \cK u  - \sigma_\fa  u\big) \rd\bx \rd\bv-\int_{\Gamma_-} \bv\cdot \bn G^\fa_b f \rd S(\bx) \rd\bv.
\end{equation}

Let us now introduce the adjoint volume Green function $G^\fa_v(\bx,\bv;\hat\bx,\hat\bv)$ as the solution of the following adjoint free transport equation:
\begin{equation}\label{EQ:Green A V}
  \begin{array}{rcll}
  	-\bv\cdot\nabla G^\fa_v(\bx,\bv;\hat\bx,\hat\bv) & =& \delta(\bx-\hat\bx)\delta(\bv-\hat{\bv}), &\mbox{in}\ \ X\\
    G^\fa_v(\bx,\bv;\hat\bx,\hat\bv) &=& 0, &\mbox{on}\ \ \Gamma_+.
  \end{array}
\end{equation}
If we multiply ~\eqref{EQ:ERT} by $G^\fa_v(\bx,\bv;\hat\bx,\hat\bv)$, subtract it by the multiplication of ~\eqref{EQ:Green A V} and $u(\bx,\bv)$, and then integrate over phase space $X$, we have
\begin{equation}\label{EQ:Sol Intgl Form B}
   u(\hat\bx,\hat\bv)= \int_X G^\fa_v \big(\sigma_\fs \cK u  - \sigma_\fa  u\big) \rd\bx \rd\bv-\int_{\Gamma_-} \bv\cdot \bn G^\fa_v f \rd S(\bx) \rd\bv,\quad (\hat\bx,\hat\bv)\in X.
\end{equation}

It turns out that we can rearrange equations~\eqref{EQ:Data A B} and ~\eqref{EQ:Sol Intgl Form B} into the form of the system~\eqref{EQ:ERT SOM A} and~\eqref{EQ:ERT SOM B} after introducing an intermediate variable. To do that, let us denote by $\tilde f$ a smooth extension of $f$ (which does exist; see for instance~\cite{Agoshkov-Book98} for justifications), and $\tilde \bn$ a smooth extension of the vector $\bn$ in the neighborhood of $\partial\Omega$. We define
\begin{equation}\label{EQ:U Def A}
  \fu(\bx,\bv)=\sigma_\fs(\bx) \cK u(\bx,\bv)-\sigma_\fa(\bx) u(\bx,\bv)-\bv\cdot\tilde\bn \tilde {f}(\bx)_{|\Gamma_-}.
\end{equation}
Then we can rewrite~\eqref{EQ:Data A B} as:
\begin{equation}\label{EQ:ERT SOM A A}
  j(\bx) = \cG^\fa_b \fu(\bx)\equiv \int_{X} G^\fa_b(\by,\bv;\bx) \fu(\by,\bv) \rd\by \rd\bv,
\end{equation}
and ~\eqref{EQ:Sol Intgl Form B} as
\begin{equation}\label{EQ:Sol Intgal Form BB}
  u(\bx,\bv)=\cG^\fa_v\fu(\bx,\bv) \equiv \int_X G^\fa_v(\hat\bx,\hat\bv;\bx,\bv) \fu(\hat\bx,\hat\bv) \rd\hat\bx \rd\hat\bv.
\end{equation}
Using~\eqref{EQ:Sol Intgal Form BB}, we can rewrite~\eqref{EQ:U Def A} as:
\begin{equation}\label{EQ:ERT SOM A B}
	\fu(\bx,\bv) = \big(\sigma_\fs \cK \cG^\fa_v - \sigma_\fa \cG^\fa_v\big) \fu(\bx,\bv) -\bv\cdot\tilde\bn \tilde {f}(\bx)_{|\Gamma_-}.
\end{equation}
Equations~\eqref{EQ:ERT SOM A A} and ~\eqref{EQ:ERT SOM A B} then form a system of exactly the same form as system~\eqref{EQ:ERT SOM A} and~\eqref{EQ:ERT SOM B}. This is the foundation of the algorithm that we will develop later. Note that the intermediate variable $\fu$ that we introduced here is a functional of both the unknown coefficient $\sigma_\fa$ to be reconstructed and the solution $u$ of the radiative transport equation with this coefficient.

\subsection{Recovering scattering coefficient}

We now assume that the absorption coefficient $\sigma_\fa(\bx)$ is known and we intend to reconstruct the scattering coefficient $\sigma_\fs(\bx)$. In this case, the adjoint boundary Green function $G^\fs_b(\bx,\bv,\by)$ that we need solves the following transport equation
\begin{equation}\label{EQ:Green S B}
  \begin{array}{rcll}
  	-\bv\cdot\nabla G^\fs_b(\bx,\bv;\by) + \sigma_\fa(\bx) G^\fs_b(\bx,\bv;\by) &=& 0, &\mbox{in}\ \ X\\
	G^\fs_b(\bx,\bv;\by) &=& \delta(\bx-\by), & \mbox{on}\ \ \Gamma_{+}
  \end{array}
\end{equation}
Following the same procedure as before, we multiply ~\eqref{EQ:ERT} by $G^\fs_b(\bx,\bv;\by)$, subtract it by the multiplication of~\eqref{EQ:Green S B} and $u(\bx,\bv)$, and then integrate over phase space $X$, we have
\begin{equation}\label{eq:subint}
   \int_{X}\big(G^\fs_b \bv\cdot\nabla u + u\bv\cdot\nabla G^\fs_b\big) \rd\bx \rd\bv = \int_{X} \sigma_\fs G^\fs_b \cK u \rd\bx \rd\bv.
\end{equation}
An integration-by-part on the left hand side then leads to
\begin{equation}\label{EQ:Data S B}
	j(\by)=\int_{X} \sigma_\fs G^\fs_b(\bx,\bv;\by) \cK u \rd\bx \rd\bv - 
	\int_{\Gamma_-} \bv\cdot \bn G^\fs_b f \rd S(\bx) \rd\bv,
\end{equation}
We now define
\begin{equation}\label{EQ:U Def S}
  \fu(\bx,\bv)=\sigma_\fs(\bx) \cK u(\bx,\bv)-\bv\cdot\tilde\bn \tilde f(\bx)_{|\Gamma_-},
\end{equation}
which then enabless us to rewrite~\eqref{EQ:Data S B} as
\begin{equation}\label{EQ:ERT SOM S A}
	j(\bx)= \cG^\fs_b \fu (\bx) \equiv \int_{X}  G^\fs_b(\by,\bv;\bx) \fu(\by,\bv) \rd\by \rd\bv.
\end{equation}

To derive the equation for $\fu$, we introduce the volume adjoint Green function $G^\fs_v(\bx,\bv;\hat\bx,\hat\bv)$ as the solution to the adjoint transport equation:
\begin{equation}\label{EQ:Green S V}
  \begin{array}{rcll}
  	-\bv\cdot\nabla G^\fs_v(\bx,\bv;\hat\bx,\hat\bv) + \sigma_\fa(\bx) G^\fs_v(\bx,\bv;\hat\bx,\hat\bv) &=& \delta(\bx-\hat\bx)\delta(\bv-\hat\bv), &\mbox{in}\ \ X\\
    G^\fs_v(\bx,\bv;\hat\bx,\hat\bv) &=& 0, &\mbox{on}\ \ \Gamma_+
  \end{array}
\end{equation}
with $(\hat\bx,\hat\bv)\in X$. This Green function allows us to derive the following result, following exactly the same procedure as before:
\begin{equation}\label{EQ:Sol Intgal Form BV}
	u(\hat\bx,\hat\bv)=\int_{X}\sigma_\fs G^\fs_v \cK u  \rd\bx \rd\bv - \int_{\Gamma_-} \bv\cdot\bn f G^\fs_v \rd S(\bx) d\bv =\int_{X} G^\fs_v \fu(\bx) \rd\bx \rd\bv \equiv \cG^\fs_v\fu(\hat\bx).
\end{equation}
Plugging ~\eqref{EQ:Sol Intgal Form BV} into ~\eqref{EQ:U Def S}, we can obtain the state equation:
\begin{equation}\label{EQ:ERT SOM S B}
	\fu(\bx,\bv) = \sigma_\fs(\bx) \cK \cG^\fs_v \fu(\bx,\bv)-\bv\cdot\tilde\bn\tilde f(\bx)_{|\Gamma_-}.
\end{equation}

Equations~\eqref{EQ:ERT SOM S A} and ~\eqref{EQ:ERT SOM S B} then form a system of exactly the same form as the system of~\eqref{EQ:ERT SOM A} and~\eqref{EQ:ERT SOM B}. As in the previous case, the intermediate variable $\fu$ that we introduced here is a functional of both the unknown coefficient $\sigma_\fs$ to be reconstructed and the solution $u$ of the radiative transport equation with this coefficient.

Let us mention that the idea of splitting of streaming (including the absorption) with the rest of the operator has been explored by Bal and Monard in~\cite{BaMo-JCP10} where they constructed an accurate solver for the streaming operator based on numerical rotations.

\subsection{Generalizations}

The formulation we presented above can be easily generalized to the case when the measured data take another form. One typical type of data assumed in the literature is given as, see for instance~\cite{BaTa-SIAM07,ChSt-IP96},
\begin{equation}\label{EQ:Data 2}
	j(\bx,\bv)=\tilde\cM u(\bx) \equiv u(\bx,\bv)\delta(\bx-\bx')\delta(\bv-\bv'),\quad (\bx',\bv')\in \Gamma_+.
\end{equation}
In this case, we follow exactly the same procedures as before to get the formulation~\eqref{EQ:ERT SOM A A} and~\eqref{EQ:ERT SOM A B}, and ~\eqref{EQ:ERT SOM S A} and~\eqref{EQ:ERT SOM S B}. The only change needed would be to replace the $\delta(\bx-\by)$ terms in the equations for the boundary adjoint Green functions $G^\fa_b$ and $G^\fs_b$, i.e ~\eqref{EQ:Green A B} and ~\eqref{EQ:Green S B}, with $\delta(\bx-\bx')\delta(\bv-\bv')$. The same idea apply to the cases where the illumination source is angularly resolved, assuming it can be constructed.

We remark, however, that when data encoded in measurements of the form~\eqref{EQ:Data 2} or in illumination of the form $f(\bx,\bv)$ are available for, the inverse problem can be less ill-posed or even well-posed; see for instance~\cite{Bal-IP09,BaTa-SIAM07,ChSt-CPDE96,ChSt-IP96,Langmore-IP08,StUh-APDE08,Tamasan-IP02,Tamasan-CM03,Wang-AIHP99} and references therein. In practice, however, one can access only a very limited number of directions. Inverse problems in these settings can still be very ill-posed and the method we propose in Section~\ref{SEC:SOM} are still useful there.

\section{Matrix Representation}
\label{SEC:Disc Form}

We now construct the discrete version of the system~\eqref{EQ:ERT SOM A} and ~\eqref{EQ:ERT SOM B}. To do that, we need to discretize the radiative transport equation. Due to the fact that the unknown $u$ is posed in phase space, we need to discretize in both spatial and angular variables. There are many existing methods to perform such a discretization, see for instance~\cite{DeVo-JCP02,DuKl-JCP02,GaZh-TTSP09,KiIs-JCP99,KiMo-JCP03,KlLa-JCP06,LeMi-Book93} and references therein. We restrict ourselves to a first-order finite-volume discrete-ordinate discretization that we proposed in our earlier work~\cite{ReAbBaHi-OL04,ReBaHi-SIAM06}. 

We assume that the spatial domain of interest, $\Omega$, is discretized into a total number of $N_\Omega$ finite volume cells that centered at $\bx_1,\bx_2, \ldots,\bx_{N_\Omega}$ respectively. The angular domain, $\bbS^{d-1}$ is discretized into $N_\bbS$ directions $\bv_1,\bv_2, \ldots,\bv_{N_\bbS}$. The discrete ordinates method approximate the integral on the sphere $\bbS^{d-1}$ with the following quadrature rule
\begin{equation}\label{EQ:DOM}
   \int_{\bbS^{d-1}} u(\bx,\bv)d\bv \approx \sum_{\ell'=1}^{N_\bbS}\eta_{\ell'} u(\bx,\bv_{\ell'}),
\end{equation}
where $\eta_{\ell'}$ is the quadrature weight associated with direction $\bv_{\ell'}$. Following the same spirit, the scattering term $\cK u(\bx,\bv)$ is approximated by:
\begin{equation}\label{EQ:DOM Scatt}
   \cK u(\bx,\bv_\ell) \approx \sum_{\ell'=1}^{N_\bbS} \eta_{\ell'}k_{\ell\ell'} u(\bx,\bv_{\ell'}) - u(\bx,\bv_\ell),
\end{equation}
where $k_{\ell\ell'}=k(\bv_\ell,\bv_{\ell'})$. The normalization conditions, which are necessary to ensure the conservation of photons, take the following forms in discrete case:
\[
   \sum_{\ell'=1}^{N_\bbS} \eta_{\ell'}=1,\qquad\mbox{and}\qquad  \sum_{\ell'=1}^{N_\bbS} \eta_{\ell'} k_{\ell \ell'}=1, \quad 1\leq \ell \leq N_\bbS.
\]

In a first-order cell-centered finite-volume discretization, we approximate the spatial integration by:
\begin{equation}\label{EQ:FVM}
   \dint_{\Omega}u(\bx,\bv)d\bx \approx \sum_{m=1}^{N_\Omega} \zeta_m u(\bx_m,\bv)
\end{equation}
where $\zeta_m$ represents the volume of the $m$-th finite volume whose center is $\bx_m$.

In the presentation below, we assume that we have $N_d$ detectors in the setup, with the location of the $d$-th detector denoted by $\bz_d$, $1\le d\le N_d$. The data we measured are collected in the data vector $J\in\bbN^{N_d\times 1}$:
\[
	J=\big(j_1,\ j_2,\ \cdots,\ j_{N_{d}-1},\ j_{N_d}\big)^\ft.
\]
where $j_d=j(\bz_d)$.

\subsection{Recovering absorption coefficient}

With the discretization method introduced above, we get the discretized form field equation ~\eqref{EQ:ERT SOM A A} as:
\begin{equation}\label{EQ:ERT SOM A A Dis}
   j(\bz_d) = \sum_{\ell=1}^{N_\bbS} \sum_{m=1}^{N_\Omega} \eta_\ell \zeta_m G^\fa_b(\bx_m,\bv_\ell;\bz_d) \fu (\bx_m,\bv_\ell) .
\end{equation}
We then collect this datum for all detectors to have the following system
\begin{equation}\label{EQ:ERT SOM A A Dis MAT}
   J = \bG^\fa_b (\bH\otimes\bS) U,
\end{equation}
where $\otimes$ represents the direct product of two matrices. The vector $U\in\bbR^{N_\bbS N_\Omega\times 1}$ contains the discrete values of $\fu(\bx,\bv)$ and is arranged as:
\[
U=(\fU_1\ \fU_2\ \ldots \ \fU_{N_\bbS})^\ft,\ \mbox{with}\ \fU_\ell =(\fu(\bx_1,\bv_\ell)\ \fu(\bx_2,\bv_\ell)\ \ldots\ \fu(\bx_{N_\Omega}, \bv_\ell) ),
\]
where the superscript $\ft$ is used to denote the transpose operation. The matrix $\bG^\fa_b\in\bbR^{N_d\times N_\bbS N_\Omega}$, coming from the Green function, takes the form
\begin{equation}\label{EQ:Green A B MAT}
\bG^\fa_b=\left (\begin{array}{cccc} G_{11} & G_{12} & \ldots & G_{1N_\bbS}\\
G_{21} & \ddots & & \vdots\\ \vdots &  & \ddots & \vdots\\
G_{N_d 1} & \ldots & \ldots & G_{N_d N_\bbS} \end{array}\right),\ \mbox{where}\ G_{d\ell}=(G_b^\fa(\bx_1,\bv_\ell;\bz_d)\ \ldots\ G_b^\fa(\bx_{N_\Omega},\bv_\ell;\bz_d)).
\end{equation}
This means that the elements of the matrix $\bG_b^\fa$ are
\[
(\bG^\fa_b)_{d,(\ell-1)N_\Omega+m}=G^\fa_b(\bx_m;\bz_d,\bv_\ell),\quad d=1,\ldots,N_d,\ \ell=1,\ldots, N_\bbS,\  m=1,\ldots,N_\Omega.
\]
The diagonal matrices $\bS\in\bbR^{N_\Omega\times N_\Omega}$ and $\bH\in\bbR^{N_\bbS\times N_\bbS}$ are defined respectively as
\begin{equation}\label{EQ:H S}
\bH=\left(\begin{array}{ccc}\eta_1 & & \\ & \ddots & \\ & & \eta_{N_\bbS}\end{array}\right),\quad 
\bS=\left(\begin{array}{ccc}\zeta_1 & & \\ & \ddots & \\ & & \zeta_{N_\Omega}\end{array}\right).
\end{equation}

Similarly, the discretized version of the state equation given by ~\eqref{EQ:ERT SOM A B} is:
\begin{equation}\label{EQ:ERT SOM A B Dis}
\begin{split}
	\fu(\bx_m,\bv_\ell) = &\sigma_\fs(\bx_m) \sum_{\ell'=1}^{N_\bbS} \eta_{\ell'}k_{\ell\ell'} \sum_{m'=1}^{N_\Omega} \sum_{\ell'' =1}^{N_\bbS} \eta_{\ell''}\zeta_{m'} G^\fa_v(\bx_{m'},\bv_{\ell''};\bx_m,\bv_{\ell'}) \fu(\bx_{m'},\bv_{\ell''})  \\
                     & - \sigma(\bx_m) \sum_{m'=1}^{N_\Omega} \sum_{\ell''=1}^{N_\bbS} \eta_{\ell''} \zeta_{m'} G^\fa_v(\bx_{m'},\bv_{\ell''};\bx_m,\bv_{\ell}) \fu(\bx_{m'},\bv_{\ell''})-\bv_\ell\cdot\tilde\bn\tilde{f}(\bx_m)_{|\Gamma_-},
\end{split}
\end{equation}
where $m=1,\ldots,N_\Omega$, $\ell=1,\ldots,N_\bbS$. We can write it in vector form as
\begin{equation}\label{EQ:ERT SOM A B Dis Vec}
\begin{split}
   U &= (\bK \otimes \bSigma_\fs) \bG^\fa_v (\bH \otimes \bS) U - (\bI_{N_\bbS} \otimes \bSigma) \bG^\fa_v (\bH \otimes \bS) U - F \\
       &= (\bK \otimes \bSigma_\fs - \bI_{N_\bbS} \otimes \bSigma) \bG^\fa_v (\bH \otimes \bS) U - F,
\end{split}
\end{equation}
where $\bK\in \bbR^{{N_\bbS} \times {N_\bbS}}$ is a matrix contains the discretized scattering kernel, with elements $(\bK)_{\ell\ell'}=k_{\ell\ell'}$ and $\bI_{N_\bbS} \in \bbR^{{N_\bbS} \times {N_\bbS}}$ denotes the identity matrix. The diagonal matrix $\bSigma\in\bbR^{N_\Omega\times N_\Omega}$ contains the values of the total absorption coefficient $\sigma=\sigma_\fa+\sigma_\fs$ at the center of the volume elements of the spatial mesh, $\bSigma=\bSigma_\fa+\bSigma_\fs$, $\bSigma_\fa\in\bbR^{N_\Omega\times N_\Omega}$) and $\bSigma_\fs\in\bbR^{N_\Omega\times N_\Omega}$) given as:
\[
	\bSigma_\fa=\left(\begin{array}{ccc}\sigma_\fa(\bx_1) & & \\ & \ddots & \\ & & \sigma_\fa(\bx_{N_\Omega})\end{array}\right), \qquad 
	\bSigma_\fs=\left(\begin{array}{ccc}\sigma_\fs(\bx_1) & & \\ & \ddots & \\ & & \sigma_\fs(\bx_{N_\Omega})\end{array}\right).
\]
The matrix $\bG^\fa_v\in\bbR^{N_\Omega N_\bbS\times N_\Omega N_\bbS}$ is organized such that its elements are given as
\[
	(\bG^\fa_v)_{(\ell-1) N_\Omega+m,(\ell'-1)N_\Omega+m'}=G^\fa_v(\bx_{m'},\bv_{\ell'};\by_m,\bv_{\ell}), \qquad 1\le m,m'\le N_\Omega, \quad 1\le \ell,\ell'\le N_\bbS.
\]

We can now write down the two algebraic equations ~\eqref{EQ:ERT SOM A A Dis MAT} and ~\eqref{EQ:ERT SOM A B Dis Vec} in the form of~\eqref{EQ:ERT SOM A} and ~\eqref{EQ:ERT SOM B} with matrices:
\begin{eqnarray}
\label{EQ:ERT SOM A A Dis MAT A} \bA &=& \bG^\fa_b (\bH\otimes\bS), \\
\label{EQ:ERT SOM A B Dis MAT B} \bB &=&(\bK \otimes \bSigma_\fs - \bI_{N_\bbS} \otimes \bSigma) \bG^\fa_v (\bH \otimes \bS).
\end{eqnarray}
Note here the matrix $\bB$ depends on both  $\bSigma$ and $\bSigma_\fs$. Thus this formulation can also be used to recover the unknown scattering coefficient if the absorption coefficient is known.

\subsection{Recovering scattering coefficient}

Using the same discretization as in the previous section, we can transform the field equation~\eqref{EQ:ERT SOM S A} into the discretized form:
\begin{equation}\label{EQ:ERT SOM S A Dis}
   j(\bz_d) = \sum_{\ell=1}^{N_\bbS} \sum_{m=1}^{N_\Omega} \eta_\ell \zeta_{m} G^\fs_b(\bx_m,\bv_\ell;\bz_d) \fu(\bx_m,\bv_\ell).
\end{equation}
We then collect this datum for all detectors to arrive at the following algebraic system:
\begin{equation}\label{EQ:ERT SOM S A Dis MAT}
   J = \bG^\fs_b (\bH\otimes\bS) U,
\end{equation}
where $\bG^\fs_b$ has the same format as $\bG^\fa_b$ defined in~\eqref{EQ:Green A B MAT}, while $\bH$ and $\bS$ are given in~\eqref{EQ:H S}.

The discretized version of the state equation~\eqref{EQ:ERT SOM S B} can now be written as, $m=1,\ldots,N_\Omega$,\ $\ell=1,\ldots,N_\bbS$:
\begin{multline}\label{EQ:ERT SOM S B Dis}
	\fu(\bx_m,\bv_\ell) = \sigma_\fs(\bx_m) \sum_{\ell'=1}^{N_\bbS} \eta_{\ell'}k_{\ell \ell'} \sum_{i'=1}^{N_\Omega} \sum_{\ell''=1}^{N_\bbS} \eta_{\ell''} G^\fs_v(\bx_{i'},\bv_{\ell''};\bx_m,\bv_{\ell'}) \fu(\bx_{i'},\bv_{\ell''}) \zeta_{i'} \\
                     -\bv_\ell\cdot\tilde\bn \tilde f(\bx_m)_{|\Gamma_-}.
\end{multline}
In matrix form, this is:
\begin{equation}\label{EQ:ERT SOM S B Dis MAT}
   U = (\bK \otimes \bSigma_\fs) \bG^\fs_v(\bH \otimes \bS) U - F,
\end{equation}
where $\bK$, $\bSigma_\fs$ are defined as before. The matrix $\bG^\fs_v$ has exactly the same structure as the matrix $\bG^\fa_v$.

The two algebraic systems~\eqref{EQ:ERT SOM S A Dis MAT} and ~\eqref{EQ:ERT SOM S B Dis MAT} are in the form of~\eqref{EQ:ERT SOM A} and ~\eqref{EQ:ERT SOM B} if we define
\begin{eqnarray}
  \label{EQ:ERT SOM S A Dis MAT A}	\bA&=&\bG^\fs_b (\bH\otimes\bS),\\
  \label{EQ:ERT SOM S B Dis MAT B}	\bB&=&(\bK \otimes \bSigma_\fs) \bG^\fs_v (\bH \otimes \bS).
\end{eqnarray}
Note that the matrix $\bB$ in this case depends only on the scattering coefficient ($\Sigma_\fs$, i.e., the discretized form of $\sigma_\fs$) that we are interested in reconstructing. It does not depend on the absorption coefficient. This is different from the matrix $\bB$ we introduced in~\eqref{EQ:ERT SOM A B Dis MAT B}.

We remark finally that even though the discretized formulation in this section is based on a first-order finite-volume discrete-ordinate method for the radiative transport equation, the reconstruction method we introduce in next section is not limited to this discretization.

\section{Subspace-based Optimization Algorithms}
\label{SEC:SOM}

Let us recall that to solve the inverse problem of reconstructing the optical parameters, the data we are given are encoded in the map $\Lambda_{\sigma_\fa,\sigma_\fs}$. This means that for each illumination source $f$, we have the corresponding current data $\Lambda_{\sigma_\fa,\sigma_\fs} f$. In practice, we have a finite number, say $N_q$, of sources to use. The total data available to us are thus:
\begin{equation}\label{EQ:Data Multi S}
	\left\{\ f_q,\ \{\ \cM u_q(\bz_d)\ \}_{d=1}^{N_d}\ \right\}_{q=1}^{N_q},
\end{equation}
where $\square_q$ denote the quantity $\square$ associated with source $q$ ($1\le q\le N_q$). The construction in the previous sections can thus be conducted for data collected from each illumination source. We have the following system of equations:
\begin{eqnarray}
	\label{EQ:ERT SOM A q} J_q &=& \bA U_q, \\
	\label{EQ:ERT SOM B q} U_q &=& \bB U_q - F_q.
\end{eqnarray}
It is critical to realize that both matrix $\bA$ and matrix $\bB$ are independent of the illuminations (i.e. independent of the index $q$). Otherwise the computation will be very expensive as we will need to compute the SVD of $\bA$.

\subsection{Singular value decomposition}
\label{SUBSEC:SVD}

Let us first briefly recall the singular value decomposition of a non-symmetric matrix $\bA\in\bbR^{N_d\times N_\Omega N_\bbS}$. We denote by $\{\mu_d\}_{d=1}^{N_d}$ the singular values, arranged in nonincreasing order $\mu_1 \geq \mu_2 \geq \cdots \geq \mu_d \geq \ldots \geq \mu_{N_d}$, where it is assumed that there are a total number of $L$ large singular values. We denote by $\{\psi_d \in \bbR^{N_d \times 1}\}_{d=1}^{N_d}$ the left singular vectors and $\{\phi_i \in \bbR^{N_\Omega N_\bbS \times 1}\}_{i=1}^{N_\Omega N_{\bbS}}$ the right singular vectors. These vectors satisfy
\begin{equation}\label{EQ:A SVs}
	\bA \phi_d = \mu_d \psi_d,  \qquad \bA^\ft \psi_d= \mu_d \phi_d, \qquad 1\le d\le N_d .
\end{equation}
Here again the superscript $\ft$ is used to denote the transpose operation. We assume that the singular vectors are all normalized to have Euclidean norm $1$ so that $\{\psi_d\}_{d=1}^{N_d}$ and $\{\phi_i\}_{i=1}^{N_\Omega N_{\bbS}}$ are orthonormal bases for $\bbR^{N_d \times 1}$ and $\bbR^{N_\Omega N_\bbS \times 1}$ respectively. The singular value decomposition of $\bA$ can then be represented as:
\begin{equation}\label{EQ:A SVD}
  \bA = \bPsi \bLambda \bPhi^\ft.
\end{equation}
The matrix $\bPsi=[\psi_1,\ldots,\psi_{N_d}]$ consists of the left singular vectors and the matrix $\bPhi=[\phi_1,\ldots,\phi_{N_\Omega N_\bbS}]$ consists of the right singular vectors. The diagonal of the rectangular diagonal matrix $\bLambda \in \bbR^{N_d \times N_\Omega N_\bbS}$ contain the singular values. 

\subsection{Signal and noise subspaces}
\label{SUBSEC:Subspace}

Due to the instability of inverse transport problems, the matrix $\bA$ defined in~\eqref{EQ:ERT SOM A q} has many small (and zero when $\bA$ is not invertible) singular values. High-frequency noise in the data $J_q$ can easily ruin the reconstruction of $U_q$ ($1\le q\le N_q$) through these small singular values. The resolution in the reconstruction of the coefficients is thus limited by the smallest singular value, which we denote by $\mu_c$, that is not significantly ruined by noise. We assume that there are $L$ large (i.e. $\ge \mu_c$) singular values, i.e. $\mu_1 \geq \mu_2 \geq \cdots \geq \mu_L = \mu_c \gg \mu_{L+1} \geq \ldots \geq \mu_{N_d}$. 

We decompose the matrix $\bPhi$ into $\bPhi=[\bPhi^s, \bPhi^n]$ where $\bPhi^s$ and $\bPhi^n$ are matrices that contain the first $L$ and the last $N_\Omega N_\bbS-L$ columns of $\bPhi$ respectively. Then $\bPhi^s$ and $\bPhi^n$ decompose the column space of $\bPhi$ into two subspaces: the signal subspace and the noise subspace, following the terminologies in~\cite{Chen-JOSA09,PaChZhYe-JOSA10,ZhCh-IP09}. The signal subspace $V=\mbox{span}\{\phi_i\}_{i=1}^L$ is spanned by the column vectors of $\bPhi^s$, and the noise subspace $W=\mbox{span}\{\phi_i\}_{i=L+1}^{N_\Omega N_\bbS}$ is spanned by the columns of $\bPhi^n$. By construction, $W$ is the orthogonal complement of $V$ and vice versa.

Now for any known vector $U_q$ as a solution to the system~\eqref{EQ:ERT SOM A q} and ~\eqref{EQ:ERT SOM B q}, we can decompose it into a summation of its projection to the signal subspace, $U_q^s$, and its projection to the noise subspace, $U_q^n$:
\begin{equation}\label{EQ:U Us Un}
  U_q = U_q^s + U_q^n =\bPhi^s \beta^s_q+\bPhi^n \beta^n_q,
\end{equation}
where $\beta^s_q=(\beta_{q,1},\beta_{q,2},\cdots,\beta_{q,L})^\ft$ and $\beta_q^n=(\beta_{q,L+1},\beta_{q,L+2},\cdots,\beta_{q,N_\Omega N_\bbS})^\ft$ are the corresponding coefficient vectors for the projections, with $\beta_{q,i}=U_q^\ft\phi_i$ ($1\le i\le N_\Omega N_\bbS$).

\subsection{A two-step subspace minimization algorithm}
\label{SUBSEC:TSSOM}

To solve the inverse transport problems to find the optical coefficients, we first need to propagate information contained in the current data $J_q$ $(1\le q\le N_q)$ to the intermediate quantity, $U_q$, by ``solving'' the equation~\eqref{EQ:ERT SOM A q}. Due to the the existence of small singular values ($<\mu_c$) for the discrete operator $\bA$, the high-frequency components of $U_q$, i.e., the components in the noise subspace, can not be stably reconstructed in the inversion process without additional \emph{a priori} information. We should thus focus on the reconstruction of the stable components, the components in the signal subspace. 

The first algorithm we propose here completely neglects the high-frequency contents of the unknown $U_q$. In other words, we assume that:
\begin{equation}
	U_q=U_q^s=\bPhi^s\beta_q^s.
\end{equation}
This assumption allows us to construct the following two-step reconstruction process. In the first step, we reconstruct, $U_q^s$ (and thus $U_q$), from~\eqref{EQ:ERT SOM A q}. This is done by solving the minimization problem:
\begin{equation}\label{EQ:Min A}
  \{\widetilde\beta_q^s\}_{q=1}^{N_q}=\argmin_{\{\beta_q^s\in\bbR^L\}_{q=1}^{N_q}} \cO_A(\{\beta_q^s\}_{q=1}^{N_q})\equiv \sum_{q=1}^{N_q}\frac{\| (\bA\bPhi^s\beta_q^s-J_q\|_{l^2}^2}{\| J_q \|_{l^2}^2}.
\end{equation}
The solution of the minimization problem is given analytically as:
\begin{equation}\label{EQ:Usig}
  \widetilde\beta_q^s=(\frac{\psi^\ft_{1} J_q}{\mu_1},\cdots,\frac{\psi^\ft_{L} J_q}{\mu_L})^\ft, \quad\mbox{which leads to},\quad 
\wU_q^s \equiv \bPhi^s \widetilde\beta_q^s= \sum_{i=1}^{L} \frac{\psi^\ft_{i} J_q }{\mu_i} \phi_i.
\end{equation}

In the second step, we reconstruct the unknown optical coefficient through~\eqref{EQ:ERT SOM B q} using the reconstructed $U_q=\wU_q^s$ in ~\eqref{EQ:Usig}. This is done by solving the following minimization problem:
\begin{equation}\label{EQ:Min B}
 \widetilde\Sigma_\fx=\argmin_{\Sigma_\fx\in\bbR^{N_\Omega}} \cO_B(\Sigma_\fx)\equiv \sum_{q=1}^{N_q}\frac{\| (\bB-\bI_{N_\Omega N_\bbS}) \wU_q^s - F_q\|_{l^2}^2}{\| \wU_q^s \|_{l^2}^2},
\end{equation}
where the weighting factor for source $q$, $\|\wU_q^s\|_{l^2}^2=\|\widetilde\beta_q^s\|_{l^2}^2$. The unknown variable in the optimization, $\Sigma_\fx=\Sigma_\fa$ when the absorption coefficient is to be reconstructed, in which case $\bB$ is given in~\eqref{EQ:ERT SOM A B Dis MAT B}, $\Sigma_\fx=\Sigma_\fs$ when the scattering coefficient is to be reconstructed, in which case $\bB$ is given in~\eqref{EQ:ERT SOM S B Dis MAT B}. In both cases, $\bB$ is linearly related to $\Sigma_\fx$, so the minimization problem~\eqref{EQ:Min B} is a quadratic problem. This minimization problem is solved with a quasi-Newton method that we will describe briefly in Section~\ref{SUBSEC:OSSOM}.

\subsection{A modified two-step subspace minimization algorithm}
\label{SUBSEC:MTSSOM}

We can modify slightly the above two-step reconstruction algorithm to take into account part of the noise component of the unknown $U_q$. In principle, we can reconstruct the noise component of $U_q$ in the same way as we reconstruct its signal component. However, this reconstruction is not stable due to the smallness of the singular values associated with the noise component. We take here different approach. We assume that
\begin{equation}\label{EQ:Assumption B}
  U_q=\wU_q^s+\sum_{i=1}^{N_d-L} \gamma_i^n \phi_{L+i}  = \wU_q^s+\bPhi^{[n]} \gamma^n,
\end{equation}
where $\wU_q^s$ is given in~\eqref{EQ:Usig}, the coefficient vector $\gamma^n=(\gamma_{1}^n,\cdots,\gamma_{N_d-L}^n)^\ft$ is to be reconstructed, and $\bPhi^{[n]}$ is the matrix that contains the first $N_d-L$ columns of $\bPhi^n$. 

There are two characters in our assumption~\eqref{EQ:Assumption B} that make it very different from the previous studies in subspace-based reconstruction methods~\cite{Chen-JOSA09,PaChZhYe-JOSA10,ZhCh-IP09}. First, we do not include the last $N_\Omega N_\bbS-N_d$ terms of the noise component in this assumption, since $\bA\phi_k=0$, $\forall k\ge N_d+1$. Second, the coefficient vector $\gamma^n$ is \emph{independent} of the source index $q$. In other words, we look for a common ``averaged'' noise component for all $U_q$ ($1\le q\le N_q$).

In the first step of the algorithm, we reconstruct the coefficient $\gamma^n$ by solving the following minimization problem:
\begin{equation}\label{EQ:Min A Modi}
  \widetilde\gamma^n=\argmin_{\gamma^n\in\bbR^{N_d-L}} \cO_{A'}(\gamma^n)\equiv \sum_{q=1}^{N_q}\frac{\| (\bA\bPhi^{[n]} \gamma^n+\bA\wU_q^s-J_q\|_{l^2}^2}{\| J_q \|_{l^2}^2}.
\end{equation}
The solution to this problem can again be found analytically. It is
\begin{equation}\label{EQ:Unoi}
  \widetilde\gamma_i^n=\dfrac{1}{\mu_{L+i}}\sum_{q=1}^{N_q} \frac{\psi_{L+i}^\ft \widetilde J_q}{\|J_q \|_{l^2}^2}, \quad\mbox{which leads to},\quad 
\wU_q^n \equiv \bPhi^{[n]} \widetilde\gamma^n= \sum_{i=1}^{N_d-L} \big(\dfrac{1}{\mu_{L+i}}\sum_{q=1}^{N_q} \frac{\psi_{L+i}^\ft \widetilde J_q}{\|J_q \|_{l^2}^2}\big) \phi_{L+i},
\end{equation}
where $\tilde J_q=J_q-\bA\wU_q^s$. As we can see clearly, the coefficient $\widetilde\gamma_i^n$ is recovered as an average over different sources. This makes the inversion more stable.

In the second step, we reconstruct the unknown optical coefficient through~\eqref{EQ:ERT SOM B q} using the reconstructed $\widetilde U_q=\wU_q^s+\wU_q^n$ in ~\eqref{EQ:Usig}. This is done by solving the following minimization problem:
\begin{equation}\label{EQ:Min B Modi}
  \widetilde{\widetilde\Sigma}_\fx=\argmin_{\Sigma_\fx\in\bbR^{N_\Omega}} \cO_{B'}(\Sigma_\fx)\equiv \sum_{q=1}^{N_q}\frac{\| (\bB-\bI_{N_\Omega N_\bbS}) \wU_q - F_q\|_{l^2}^2}{\| \wU_q \|_{l^2}^2},
\end{equation}
where the weighting factor for source $q$, $\|\wU_q\|_{l^2}^2=\|\widetilde\beta_q^s\|_{l^2}^2+\|\widetilde\gamma^n\|_{l^2}^2$. As before, the unknown variable in the optimization, $\Sigma_\fx=\Sigma_\fa$ when the absorption coefficient is to be reconstructed, in which case $\bB$ is given in~\eqref{EQ:ERT SOM A B Dis MAT B}, $\Sigma_\fx=\Sigma_\fs$ when the scattering coefficient is to be reconstructed, in which case $\bB$ is given in~\eqref{EQ:ERT SOM S B Dis MAT B}. The minimization problem~\eqref{EQ:Min B Modi} is also a quadratic problem and is solved with a the same quasi-Newton method in Section~\ref{SUBSEC:OSSOM}.

\subsection{A one-step subspace minimization algorithm}
\label{SUBSEC:OSSOM}

In the two-step algorithms we introduced in Sections~\ref{SUBSEC:TSSOM} and ~\ref{SUBSEC:MTSSOM}, we used the first equation in the reformulated transport equation, ~\eqref{EQ:ERT SOM A q}, to determine the intermediate variable $U_q$ and then use the second equation, ~\eqref{EQ:ERT SOM B q}, to determine the unknown optical coefficient $\Sigma_\fx$. In conventional reconstruction algorithms, such as those in the references we cited, the two equations are used simultaneously to determine the unknown coefficient $\Sigma_\fx$. We now modify our two-step algorithms here to get a one-step reconstruction algorithm that is similar to those that have been developed in the literature.

We take the same assumption as in~\eqref{EQ:Assumption B}. To obtain robust reconstruction algorithms, however, we should still reconstruct the signal part of $U_q$, $U_q^s$, from the analytical expression in~\eqref{EQ:Usig}. To reconstruct the optical property $\Sigma_\fx$ and the coefficient $\gamma^n$, we solve the minimization problem:
\begin{equation}\label{EQ:Obj OSSOM}
   (\widehat\gamma^n,\widehat\Sigma_\fx)=\argmin_{(\gamma^n,\Sigma_\fx)\in\bbR^{N_d-L}\times \bbR^{N_\Omega}} \cO_{AB}(\gamma^n,\Sigma_\fx)
\end{equation}
where the objective function is essentially the summation of $\cO_A$ in~\eqref{EQ:Min A Modi} and $\cO_B$ in ~\eqref{EQ:Min B Modi}:
\begin{equation}
   \cO_{AB}(\gamma^n,\Sigma_\fx)\equiv \sum_{q=1}^{N_q}\frac{\| (\bA\bPhi^{[n]} \gamma^n-\widetilde J_q\|_{l^2}^2}{\| J_q \|_{l^2}^2}+\sum_{q=1}^{N_q}\frac{\| (\bB-\bI_{N_\Omega N_\bbS}) \bPhi^{[n]}\gamma^n - \widetilde F_q\|_{l^2}^2}{\| \wU_q^s \|_{l^2}^2}
\end{equation}
with $\widetilde J_q=J_q-\bA\wU_q^s$ and $\widetilde F_q=F_q-(\bB-\bI)\wU_q^s$. This minimization problem is more complicated to solve than its two-step correspondence in Section~\ref{SUBSEC:MTSSOM}. In practice, we can solve the two-step version and use the results as the initial guess for this one-step algorithm. The improvement we observed in the simulations we performed is not substantial.

We solve this minimization problem with a quasi-Newton scheme using the BFGS updating rule for the Hessian matrix whose initial value is set as an identity matrix. This Newton's scheme requires the gradient information of the objective function $\cO_{AB}$ with respect to the unknown $\gamma^n$ and $\Sigma_\fx$. These gradients can be easily computed as:
\begin{equation}
   \frac{\partial\cO_{AB}}{\partial\gamma_i^n} =\sum_{q=1}^{N_q}\Big[\frac{\mu_{L+i}\fA_q^*\psi_{L+i}}{\| J_q \|_{l^2}^2}+\frac{\fB_q^*(\bB-\bI_{N_\Omega N_\bbS})\phi_{L+i}}{\| \wU_q^s \|_{l^2}^2}\Big], \quad 1\le i\le N_d-L
\end{equation}
\begin{equation}
   \frac{\partial\cO_{AB}}{\partial\Sigma_{\fx,jj}} =\sum_{q=1}^{N_q}\frac{1}{\| \wU_q^s \|_{l^2}^2}\fB_q^*\frac{\partial \bB}{\partial\Sigma_{\fx,jj}}\bPhi^{[n]}\gamma^n, \quad 1\le j\le N_\Omega
\end{equation}
where $\fA_q=\bA\bPhi^{[n]} \gamma^n-\widetilde J_q$, $\fB_q=(\bB-\bI_{N_\Omega N_\bbS}) \bPhi^{[n]}\gamma^n - \widetilde F_q$, and $\Sigma_{\fx,jj}$ denotes the $j$th diagonal element of $\Sigma_\fx$. The derivative $\frac{\partial \bB}{\partial\Sigma_{\fx,jj}}$ is given by
\begin{equation}
\frac{\partial \bB}{\partial\Sigma_{\fa,jj}} =(\bK - \bI_{N_\bbS}) \otimes \fE_{jj} \bG^\fa_v (\bH \otimes \bS), \qquad\mbox{and}\qquad \frac{\partial \bB}{\partial\Sigma_{\fs,jj}} = (\bK \otimes \fE_{jj}) \bG^\fs_v (\bH \otimes \bS),
\end{equation}
respectively for the reconstruction of absorption and scattering coefficients. Here $\fE_{jj}\in\bbR^{N_\Omega \times N_\Omega}$ is a matrix whose $j$th diagonal element is $1$ but every other element is $0$.

\subsection{Implementation issues}
\label{SUBSEC:Impl}

It is easy to see from the presentation in Sections~\ref{SUBSEC:TSSOM},~\ref{SUBSEC:MTSSOM} and~\ref{SUBSEC:OSSOM} that after the SVD of the operator $\bA$ is computed, the algorithms we proposed are very fast. There is no need to solve forward and adjoint equations in the minimization process, which is very different from traditional minimization-based reconstruction algorithms. The price we pay is of course the calculation of the SVD of $\bA$, which is very expensive computationally. Fortunately, because $\bA$ does not depend on the unknowns to be reconstructed, its SVD can be precomputed off line. Moreover, we do not need to store the whole matrix $\bA$ but only the first $N_d$ left and right singular vectors. 

In practical applications of optical imaging, it is often the case that the object to be imaged, a piece of tissue or a small animal for instance, is placed inside a measurement device of regular shape, a cylinder or a cube for instance~\cite{LaFoGiDwHi-JBO07}. On can thus use a fixed discretization scheme for $\cA$ for a fixed measurement device. For instance, for the problem of reconstructing absorption coefficient described in Section~\ref{SUBSEC:Absorption}, $\bA$ is fixed once the device is fixed. Only one SVD needs to be done in the lifetime of the measurement device.

The algorithms in Sections~\ref{SUBSEC:TSSOM},~\ref{SUBSEC:MTSSOM} and~\ref{SUBSEC:OSSOM} are with increasing complexity. In practice, we can use the two-step SOM method in Section~\ref{SUBSEC:TSSOM} to construct a good approximation of the inverse solution and use it as initial guess for the modified two-step algorithm in Section~\ref{SUBSEC:MTSSOM} whose final solution can be used as an initial guess for the one-step SOM algorithm in Section~\ref{SUBSEC:OSSOM}. Our numerical simulations show that the simplest two-step algorithm Section~\ref{SUBSEC:TSSOM} gives very good approximation to the final result.

The only tunable parameter in the algorithms is the the parameter $L$ (determined by the critical singular value $\mu_c$). In principle, $L$ should be selected so that $\mu_L\approx\mu_c$ is located at the place where a large jump of the singular value occurs. When there is no obvious jump of singular values of $\bA$, $L$ should be selected mainly according to the noise level in the data. Our principle is to look at the project of the data $J$ on $\bPsi$. The coefficient of the projection, $J^*\psi_{j}$, decays fast with $j$ until it reaches the modes where random noise dominates the true signal in the data. The turning point is where $L$ is located.

\section{Numerical Experiments}
\label{SEC:Num}

We now present some numerical simulations to demonstrate the performance of the algorithms we have developed. We focus on two-dimensional settings to simplify the computation but emphasize that the discretization carries straightforwardly to three-dimensional case; see~\cite{ReBaHi-SIAM06} for typical three-dimensional reconstructions from this discretization but with a different minimization algorithm. We also non-dimensionalize the transport equation with the typical length scale of the domain and the intensity of the illumination source, so that all the numbers we show below are without unit. 

The domain we consider is the square $[0, 2]^2$. The ``measurements'' that we use in the reconstructions are synthetic data that are generated by solving the radiative transport equation for \emph{known} optical coefficients. To reduce the degree of ``inverse crimes'', we use two different sets of finite volume meshes when generating the synthetic data and when performing the numerical reconstructions. In general, the meshes for generating data are twice as fine as the meshes used in the inversion. For the noisy data, we added multiplicative random noise to the data by simply multiplying each datum by $(1+\gamma\times 10^{-2}$ \verb|random|) with \verb|random| a uniform random variable taking values in $[-1,1]$, $\gamma$ being the noise level in percentage. In each case below, we perform reconstructions using data contains noise at three different levels: (i) noiseless data ($\gamma=0$); (ii) data containing $3$\% random noise ($\gamma=3$) and (iii) data containing $10$\% random noise ($\gamma=10$). Let us emphasize that the ``noiseless'' data in (i) still contain noise that come from interpolating from the forward meshes to the inversion meshes. We use $N_q=8$ total illumination sources and for each source we measure the boundary current data at $80$ detectors uniformly distributed on $\partial\Omega$.

\begin{figure}[ht]
\centering
\includegraphics[angle=0,width=0.23\textwidth]{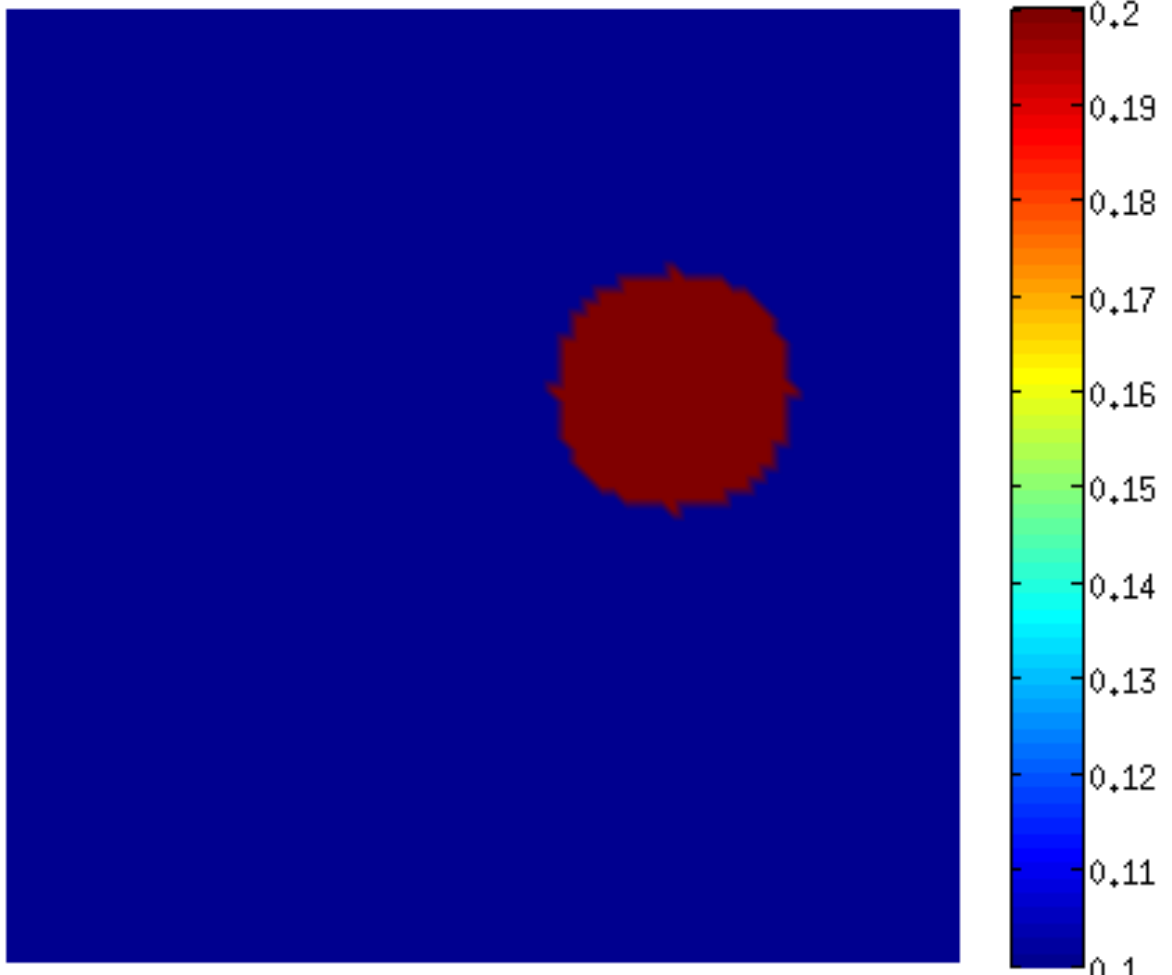}
\includegraphics[angle=0,width=0.23\textwidth]{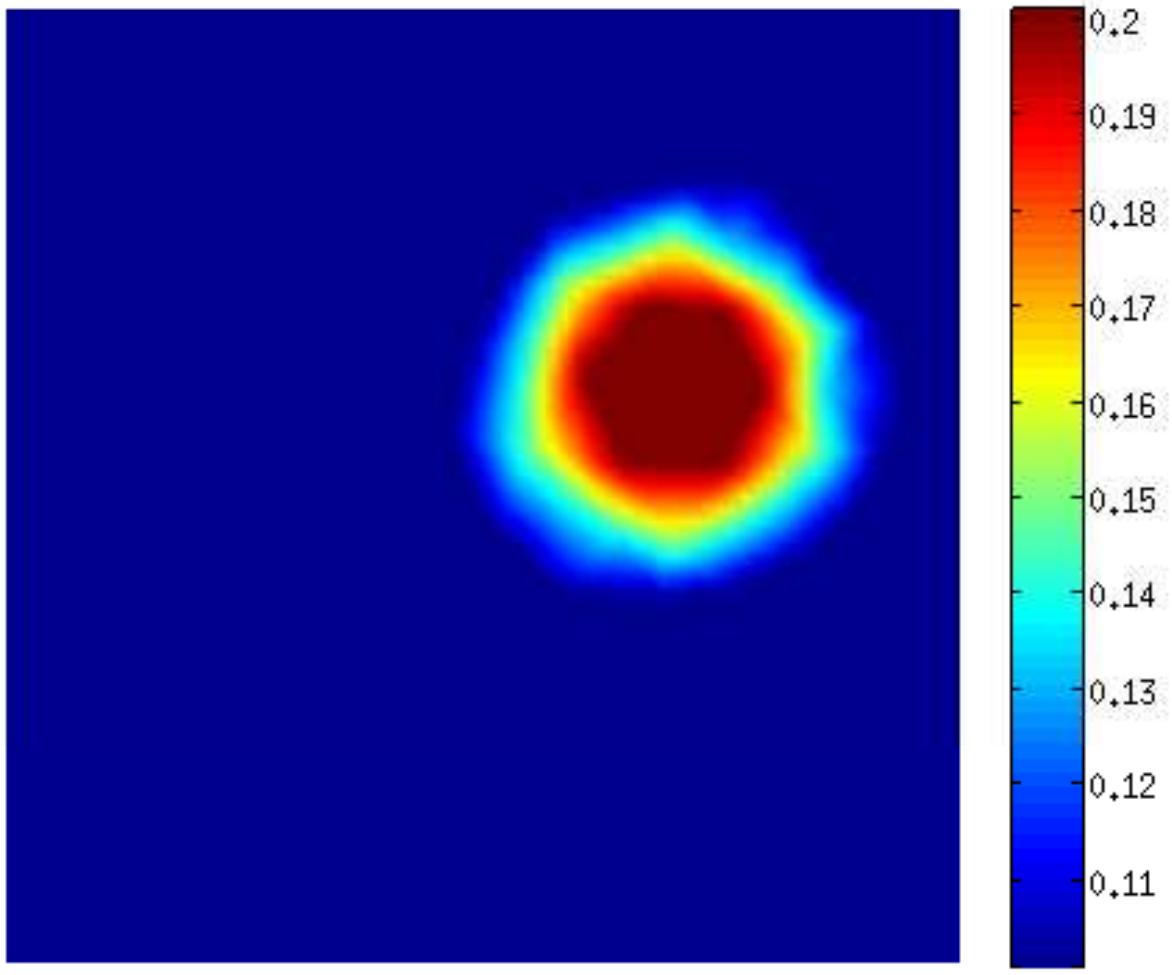}
\includegraphics[angle=0,width=0.23\textwidth]{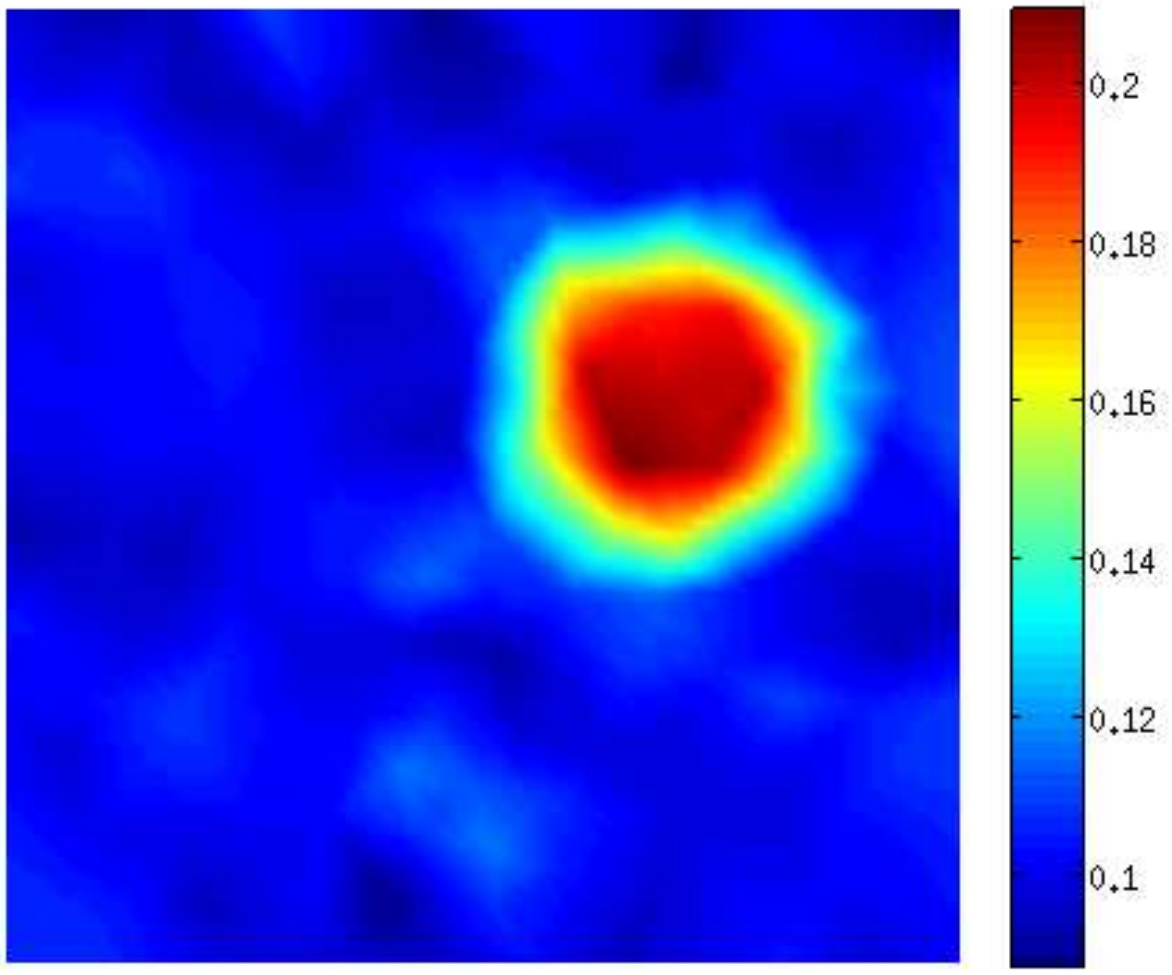}
\includegraphics[angle=0,width=0.23\textwidth]{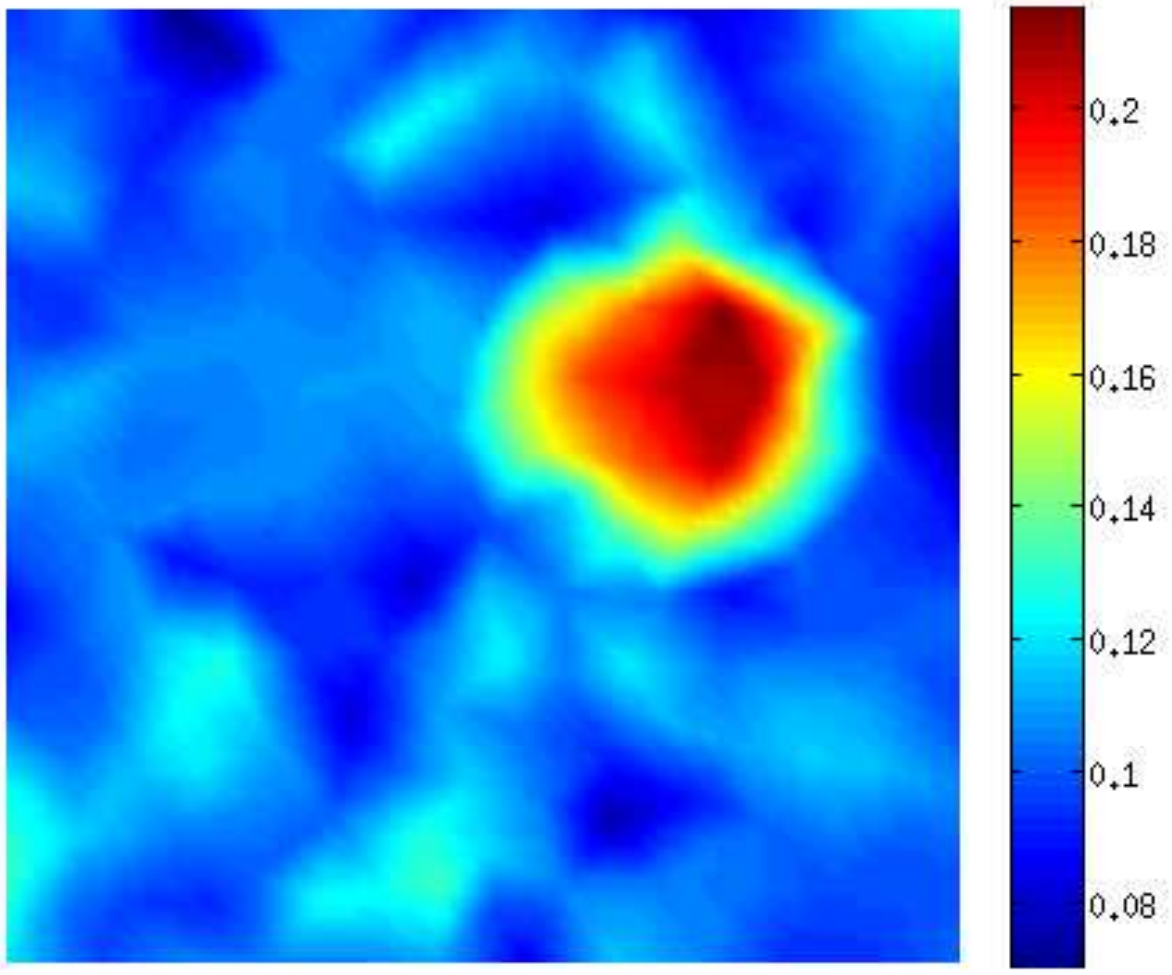}
\caption{From left to right: True absorption coefficient and absorption coefficient reconstructed with data of type (i), type (ii) and type (iii) respectively for Experiment I.}
\label{FIG:Exp Cir}
\end{figure}
\paragraph{Experiment I.} We show in Fig.~\ref{FIG:Exp Cir} the reconstructions of an absorbing disk embedded in a constant background. Shown in Fig.~\ref{FIG:Exp Cir} are the true absorption coefficient and the reconstructions with synthetic data of type (i), (ii) and (iii) respectively. The scattering coefficient is set to be constant $\sigma_s=8$. We use an isotropic scattering kernel in this case. We did not observe much difference in the reconstructions when use the Henyey-Greenstein phase function with $g=0.9$ and $\sigma_s=80$ (so that the effective scattering coefficient is still $8$). The relative $L^2$ error in the reconstructions are $2.84\%$, $5.82\%$, and $9.23\%$ for reconstructions with the three data types respectively. 

\begin{figure}[ht]
\centering
\includegraphics[angle=0,width=0.30\textwidth]{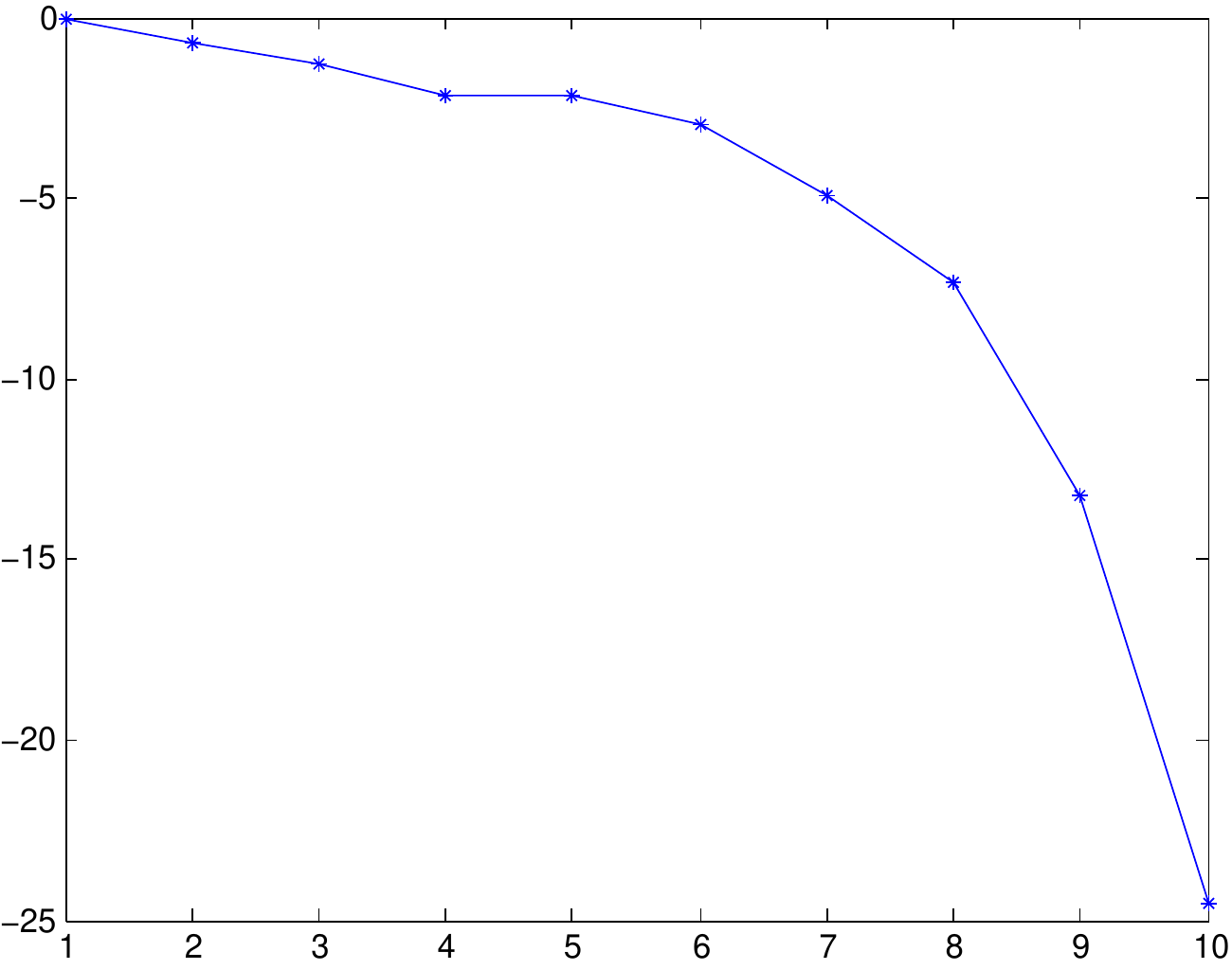}\hskip 0.6cm
\includegraphics[angle=0,width=0.30\textwidth]{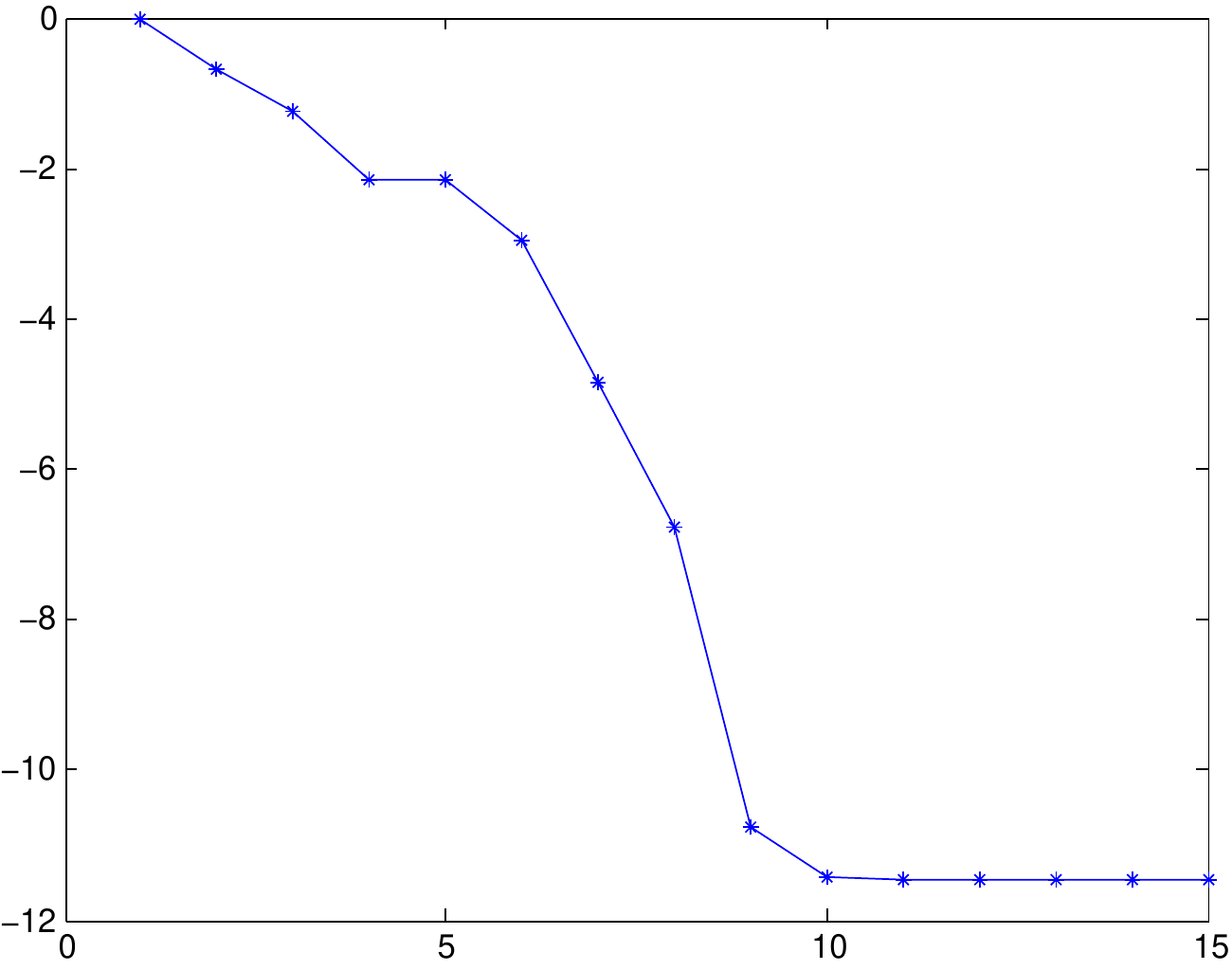}\hskip 0.6cm
\includegraphics[angle=0,width=0.30\textwidth]{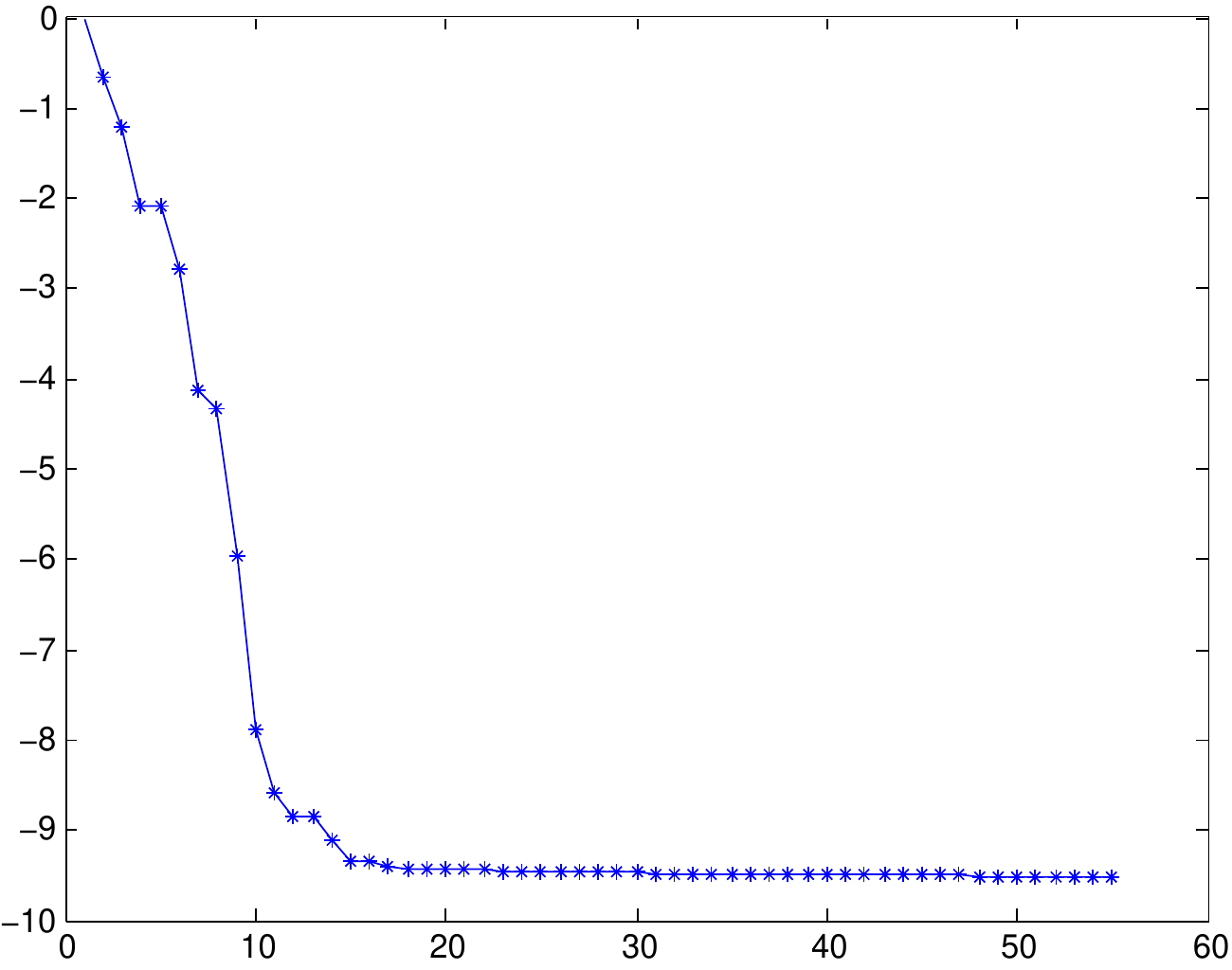}\hskip 0.6cm
\caption{Evolution of the objective functions (normalized by its starting value and in logarithmic  scale) in BFGS Newton iteration for the three reconstructions in Experiment I.}
\label{FIG:Exp Cir Con His}
\end{figure}

\begin{figure}[ht]
\centering
\includegraphics[angle=0,width=0.230\textwidth]{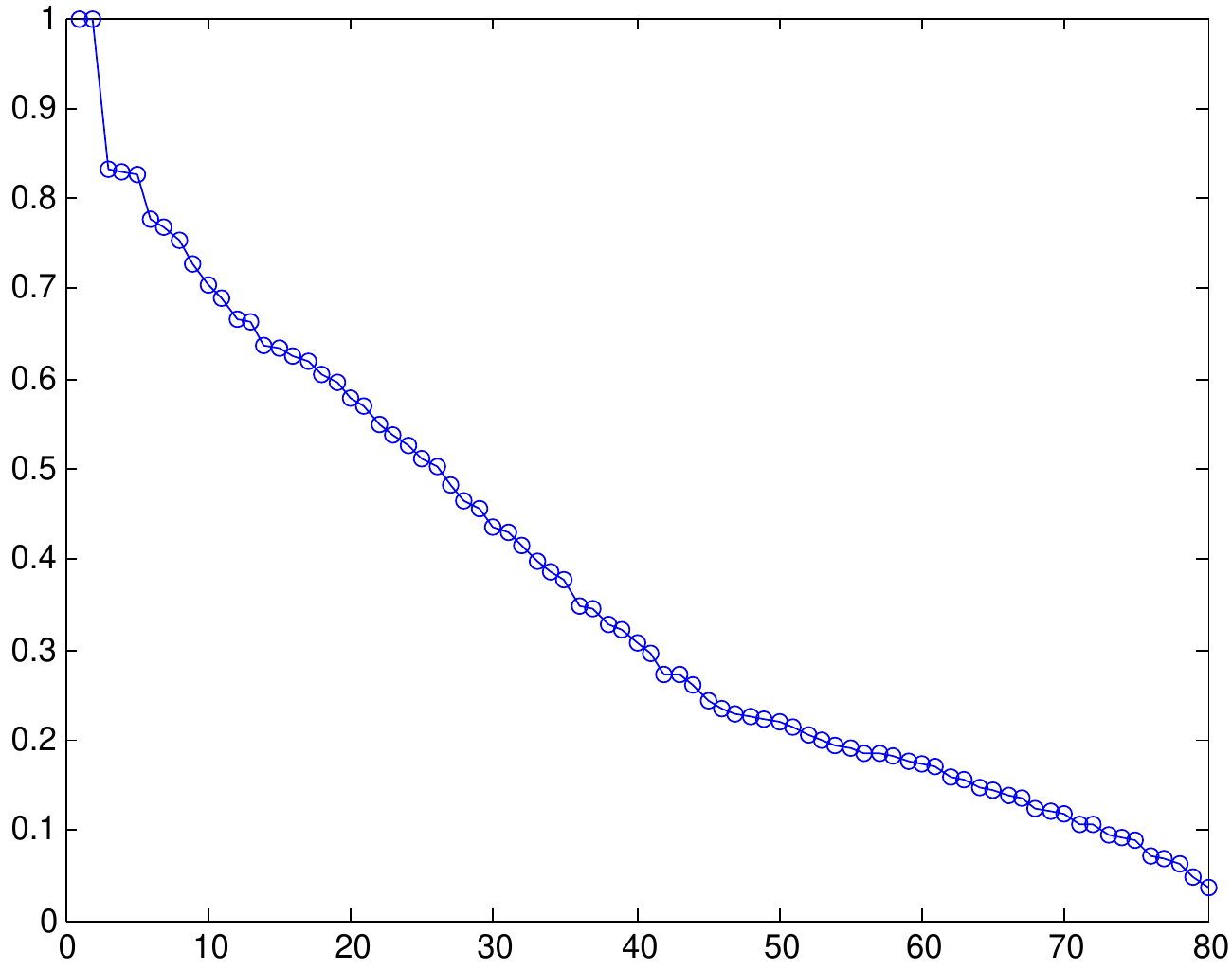}
\includegraphics[angle=0,width=0.230\textwidth]{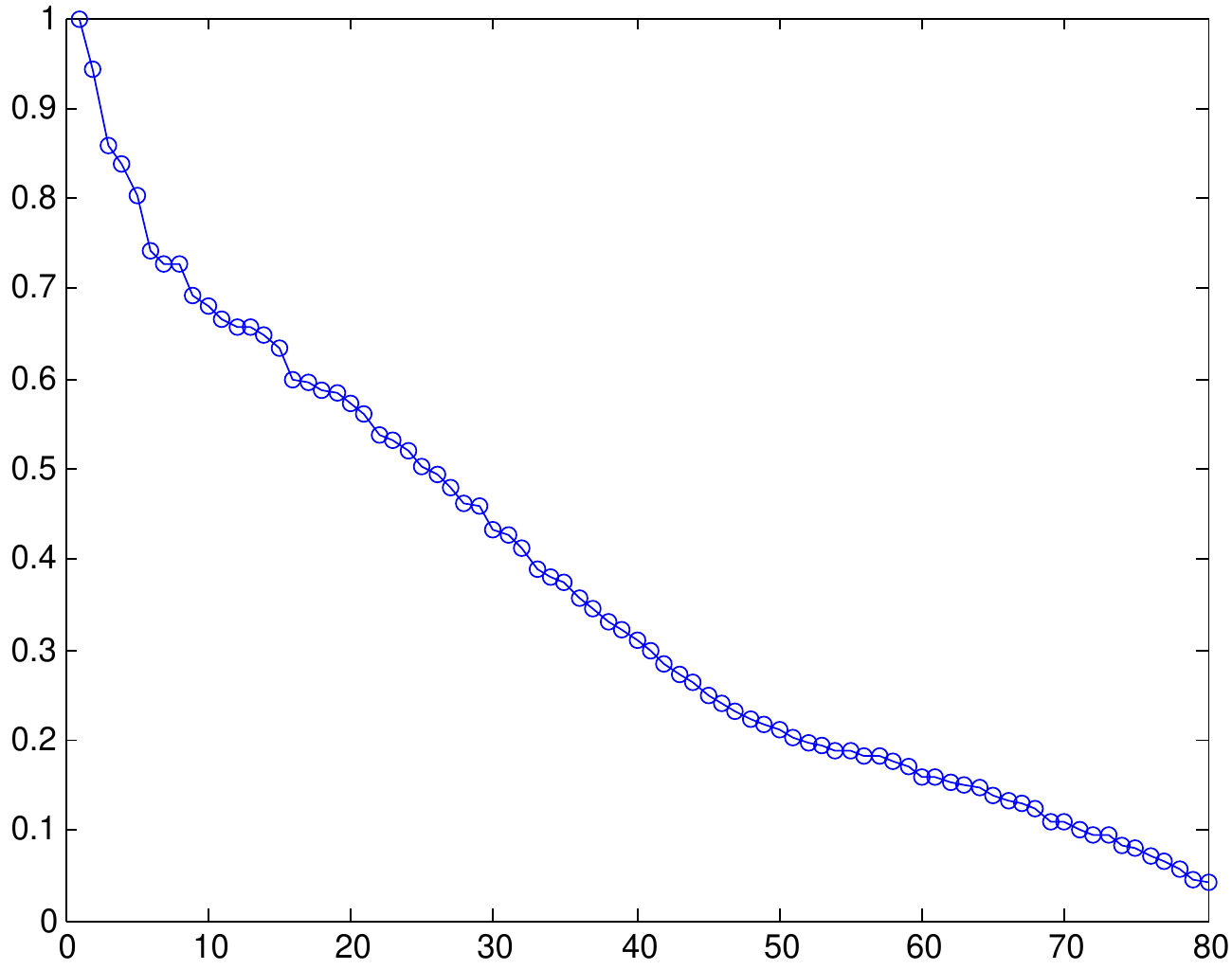}
\includegraphics[angle=0,width=0.230\textwidth]{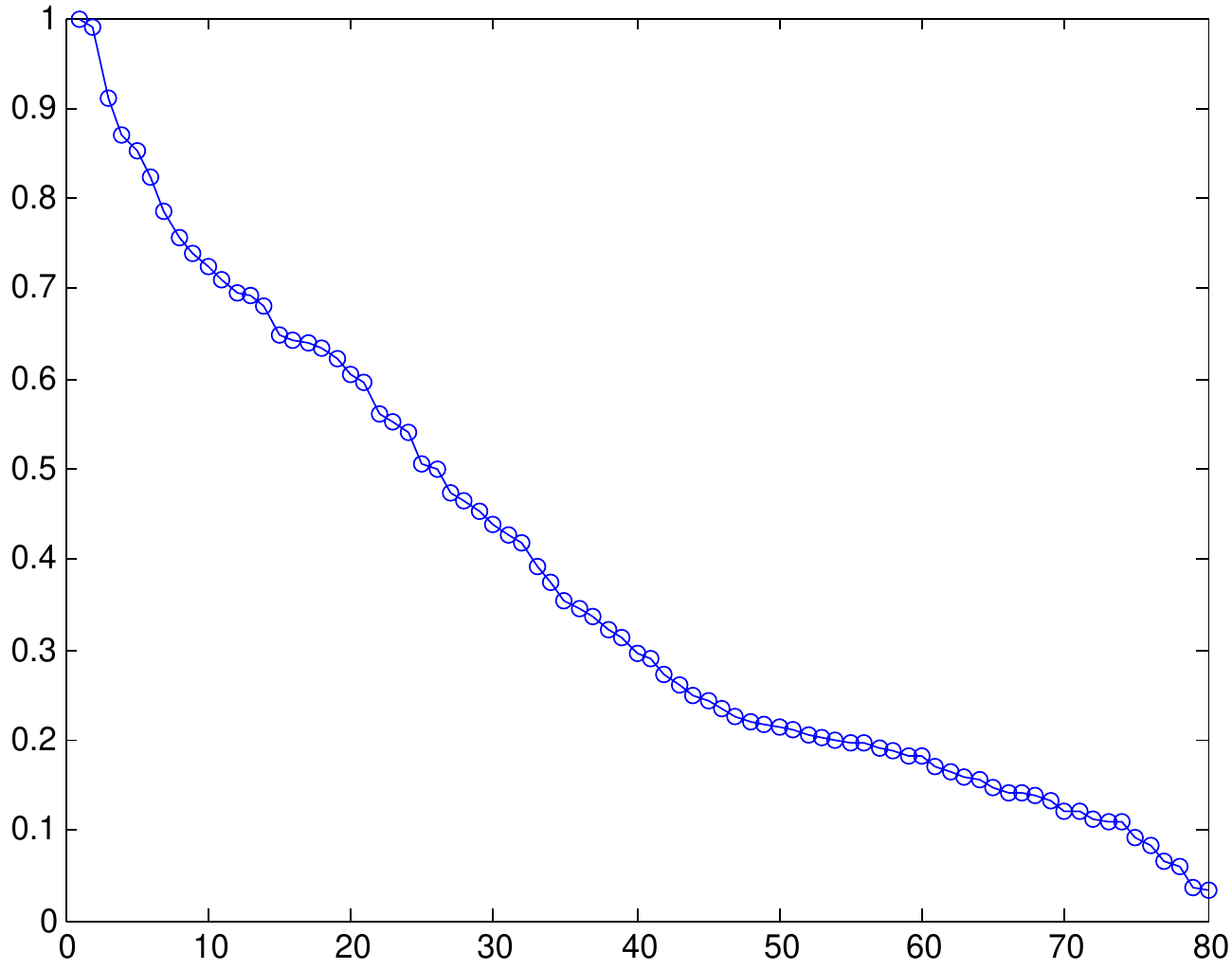}
\includegraphics[angle=0,width=0.230\textwidth]{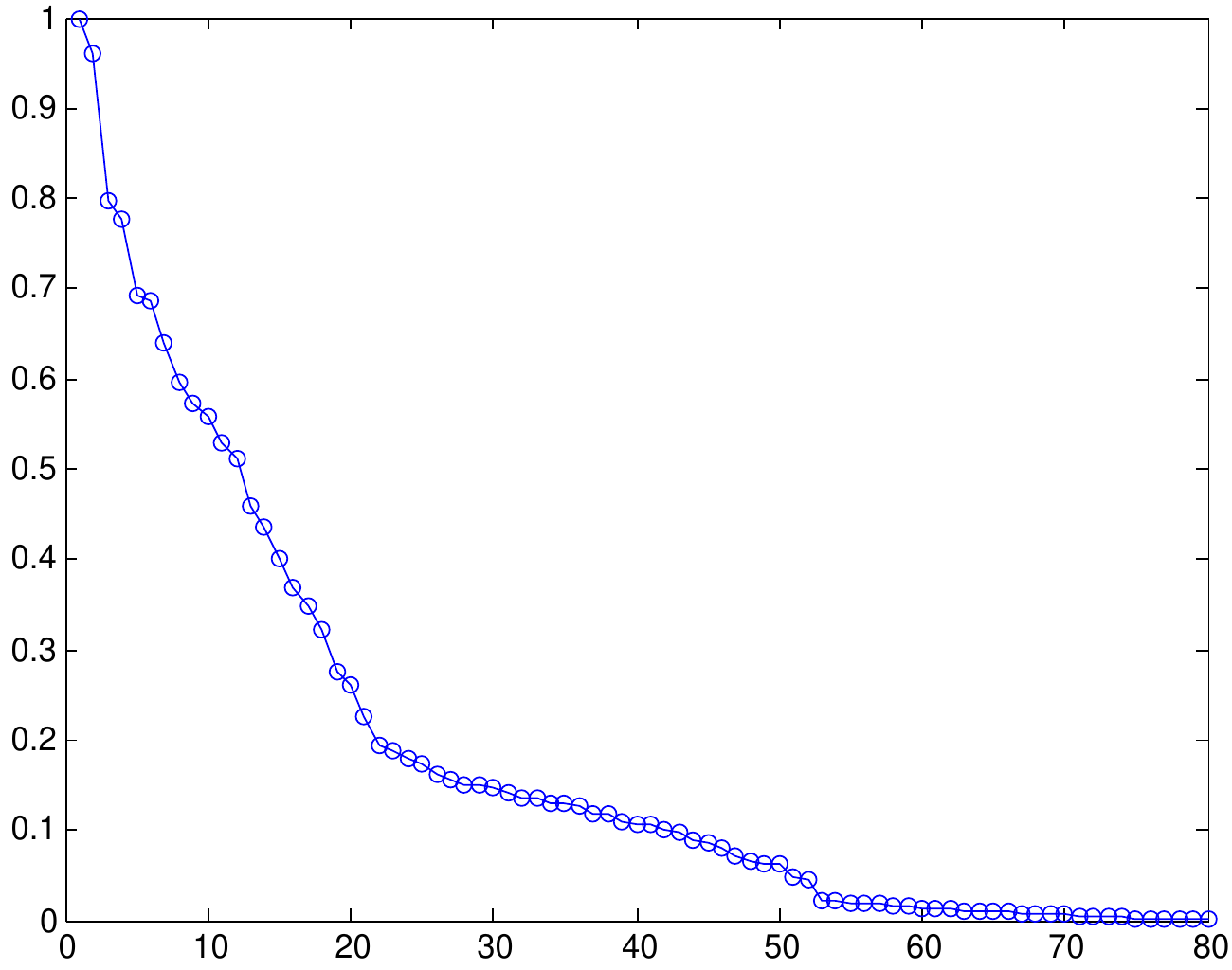}
\caption{The first $80$ singular values (normalized by the leading eigenvalue in each case) for the $\bA$ matrices in Experiment I, II, III and IV respectively.}
\label{FIG:SVD}
\end{figure}

We show in the first plot of Fig.~\ref{FIG:SVD} the singular value of the matrix $\bA$ in this case. In all the reconstructions, we take the first $L=50$ singular vectors to form the signal space and the rest to form the noise space. This works fine for the case of noiseless data. The algorithm converges in about $10$ iterations in this case as can be seen from the convergence history of the reconstruction algorithm shown in Fig.~\ref{FIG:Exp Cir Con His}. When noise is large, however, the algorithm still converges very fast at the beginning, but slow down significantly after about $10$ iterations. This happens when the algorithm struggles to find significant updates of the unknowns. Since the objective function is not getting much lower, we can stop the iteration at $10$ to get an accurate approximation to the final result. This saves significant computational time without losing much reconstruction quality. 

\begin{figure}[ht]
\centering
\includegraphics[angle=0,width=0.230\textwidth]{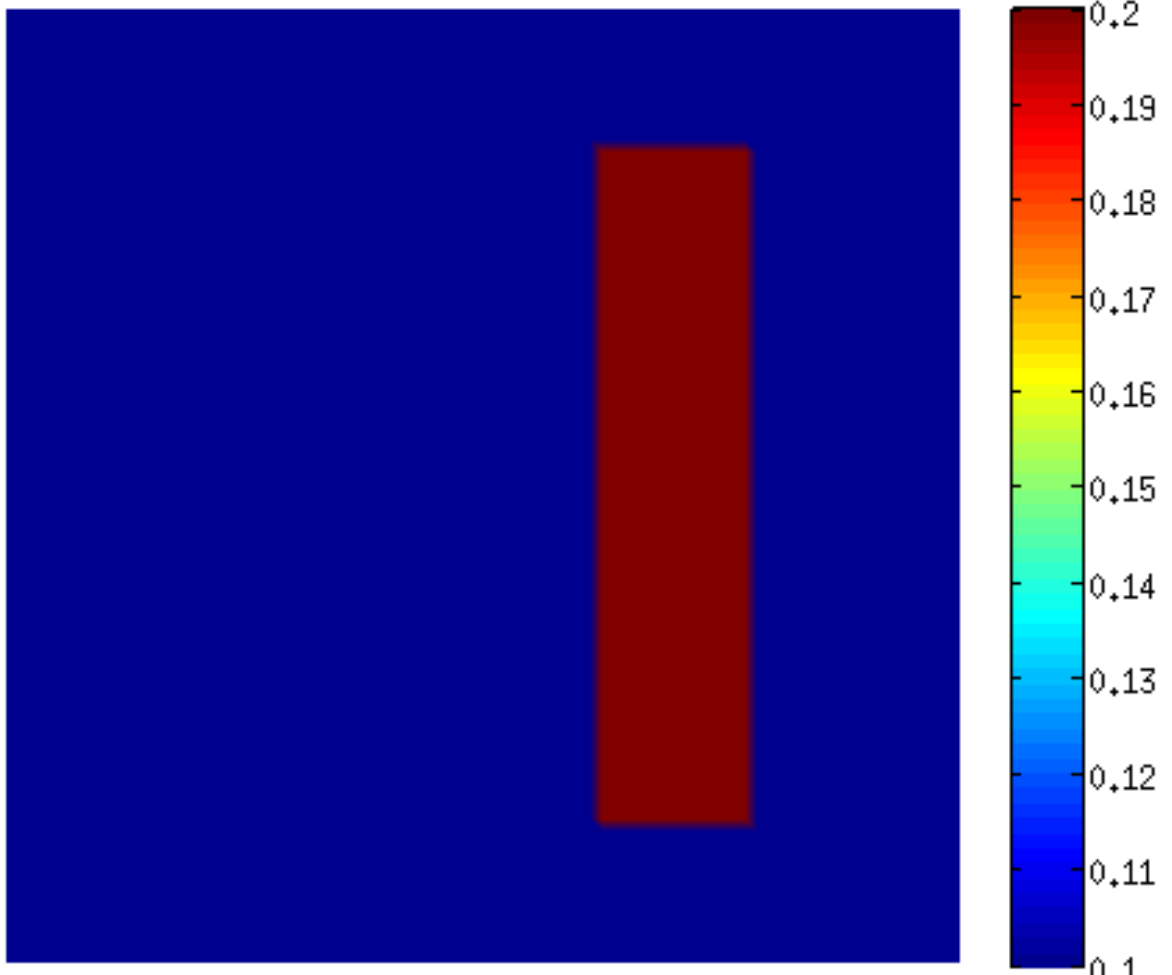}
\includegraphics[angle=0,width=0.230\textwidth]{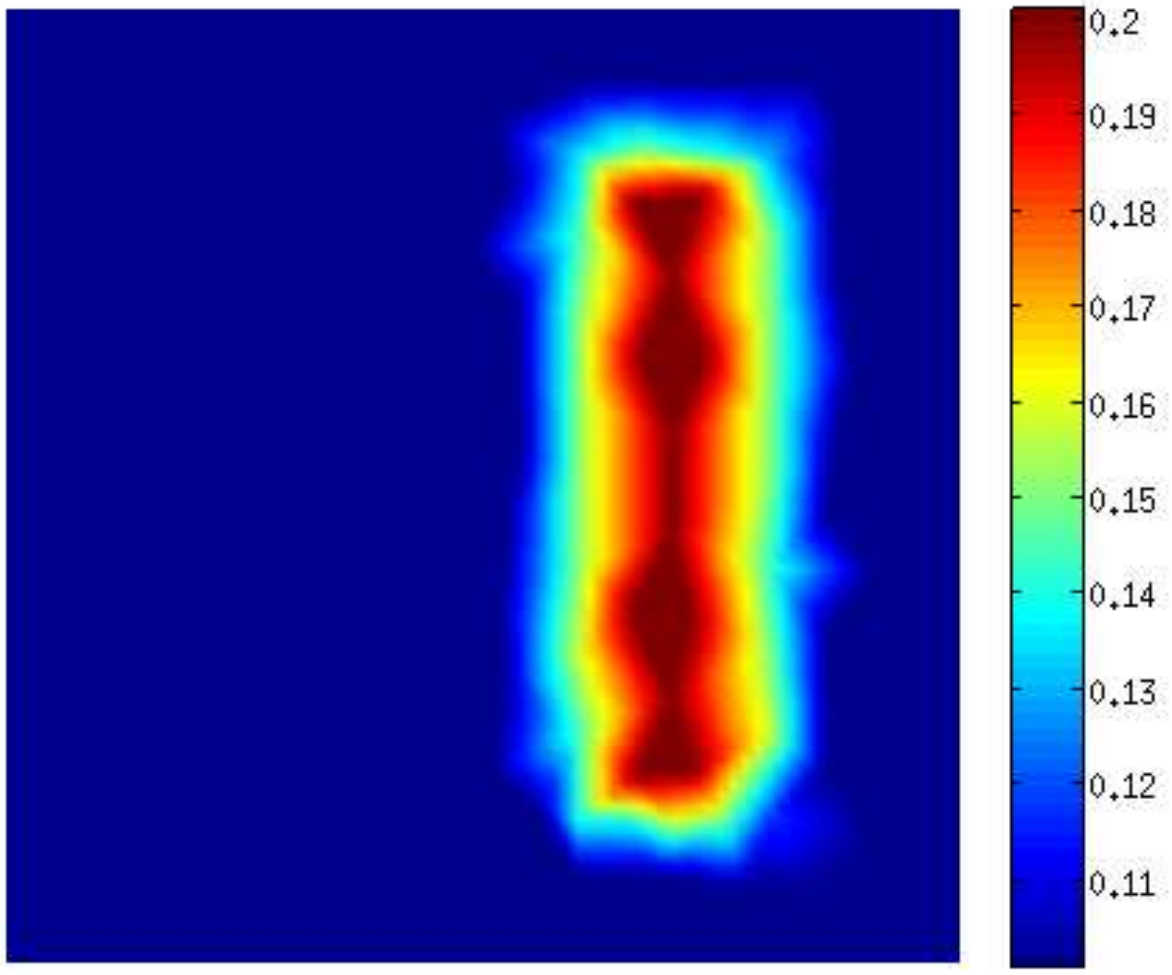}
\includegraphics[angle=0,width=0.230\textwidth]{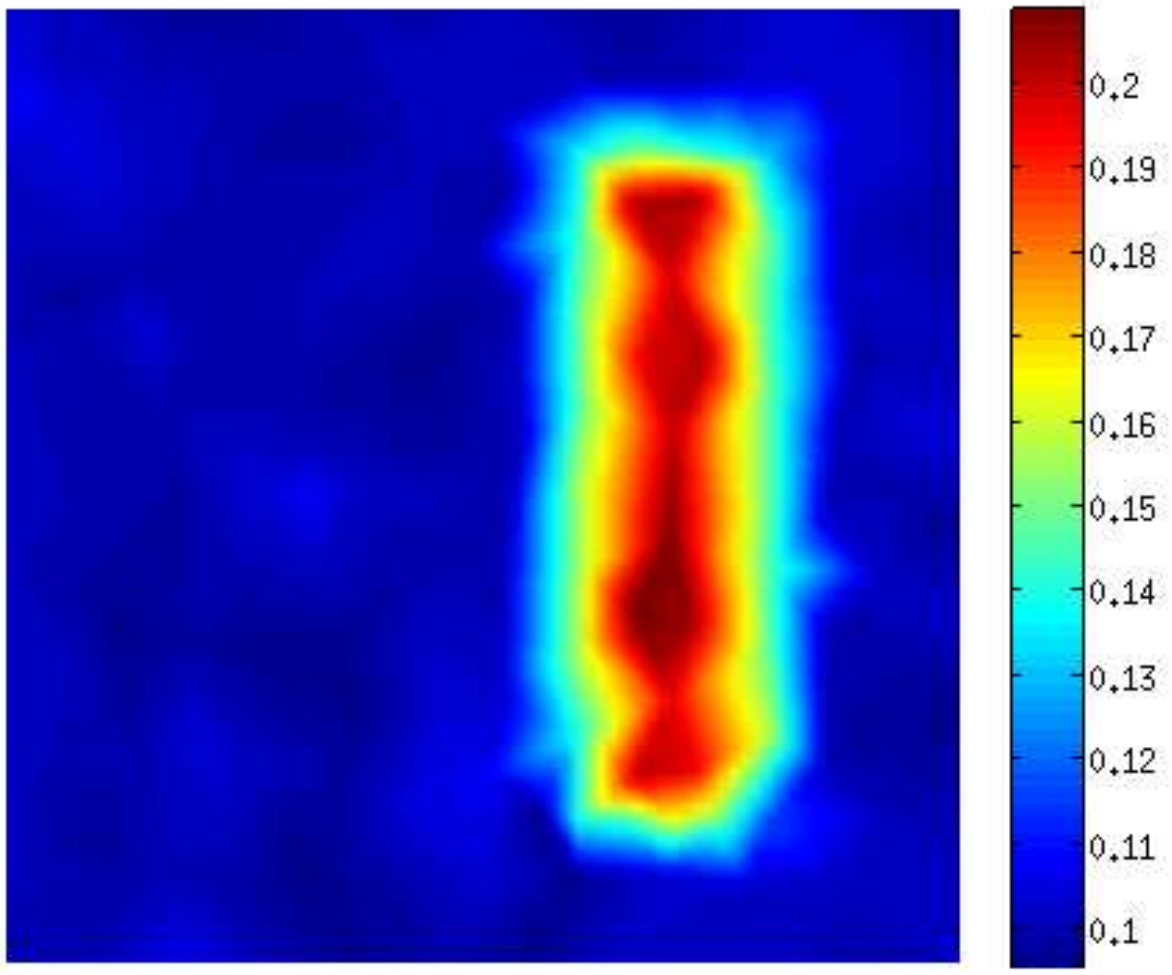}
\includegraphics[angle=0,width=0.230\textwidth]{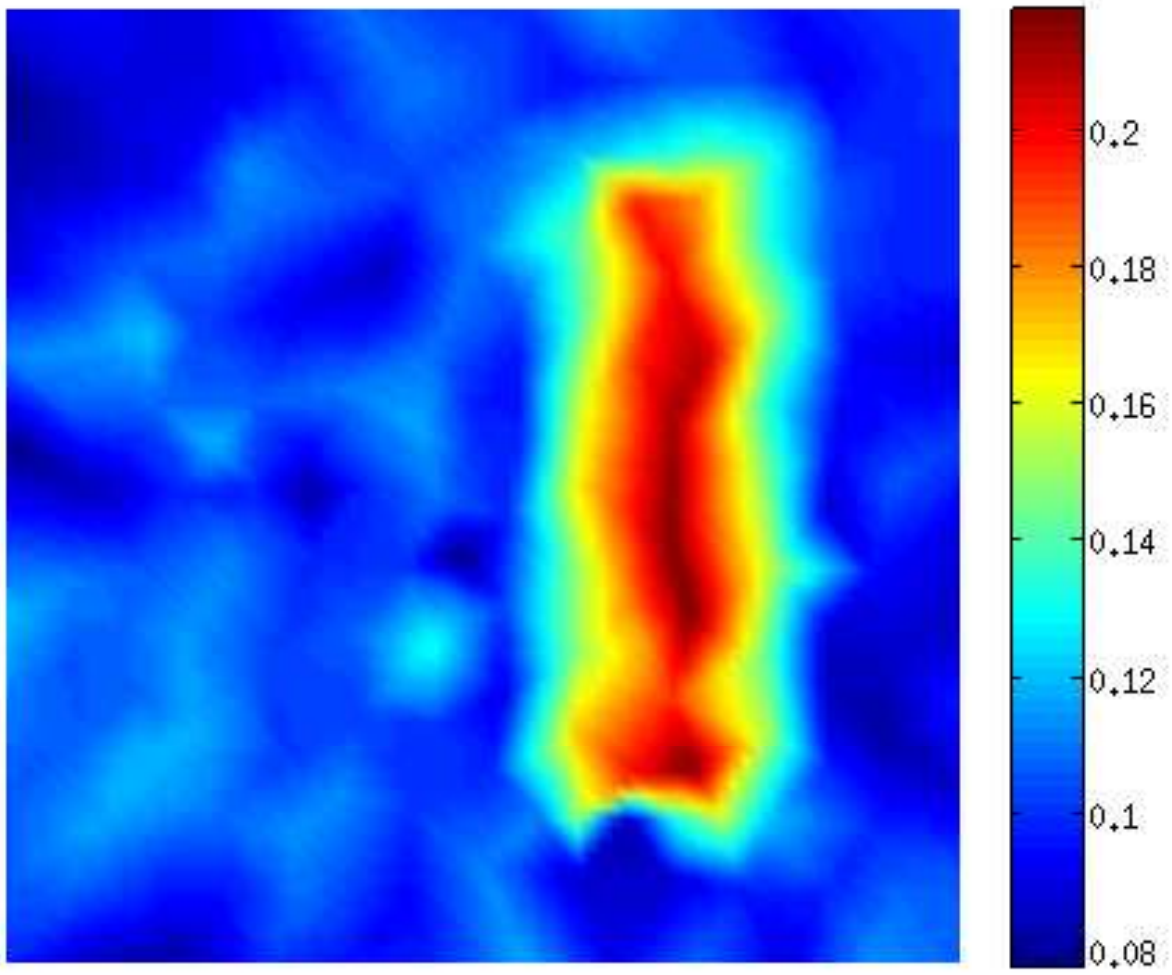}
\caption{From left to right: True absorption coefficient and absorption coefficient reconstructed with data of type (i), type (ii) and type (iii) respectively for Experiment II.}
\label{FIG:Exp Rec}
\end{figure}
\begin{figure}[ht]
\centering
\includegraphics[angle=0,width=0.30\textwidth]{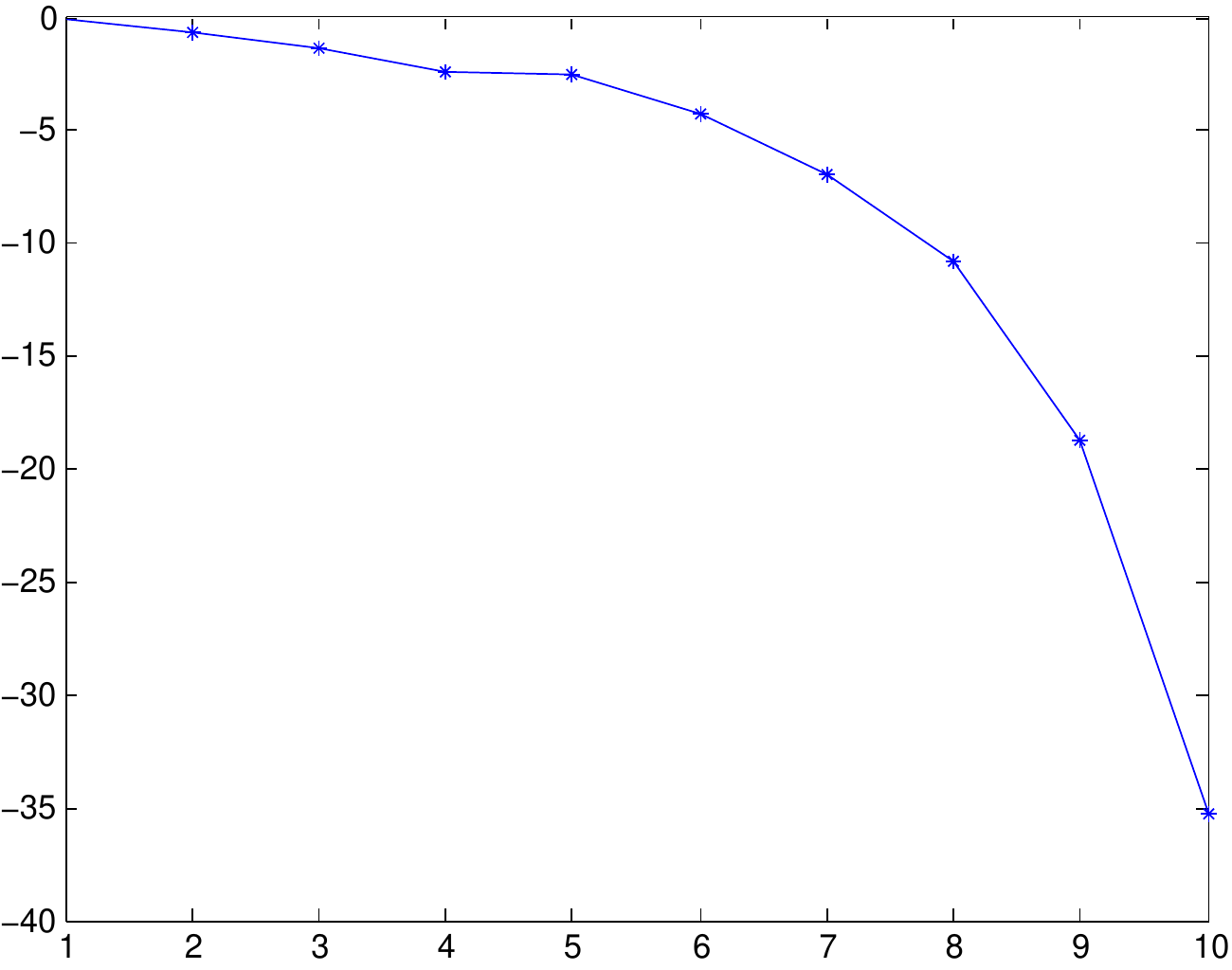}\hskip 0.6cm
\includegraphics[angle=0,width=0.30\textwidth]{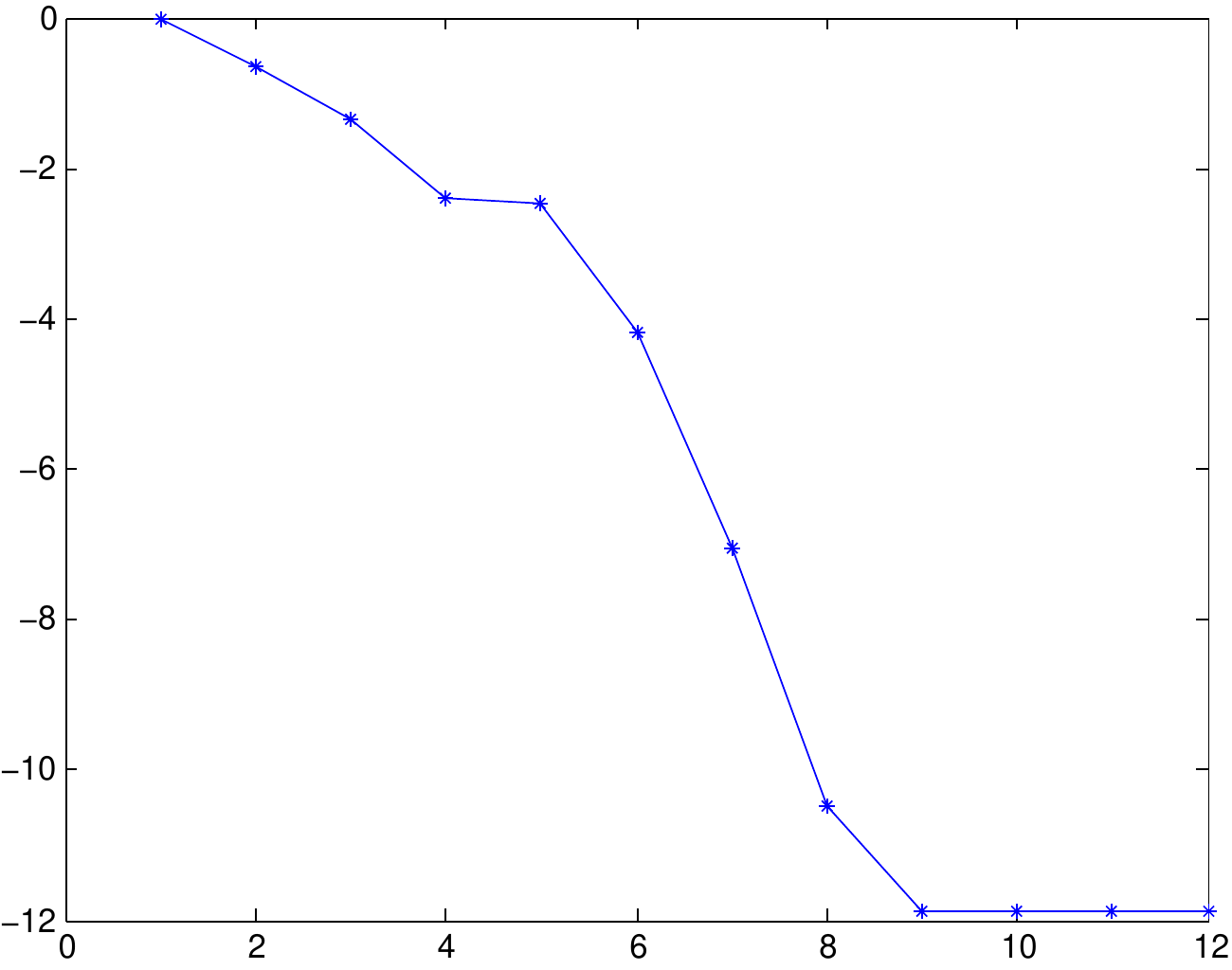}\hskip 0.6cm
\includegraphics[angle=0,width=0.30\textwidth]{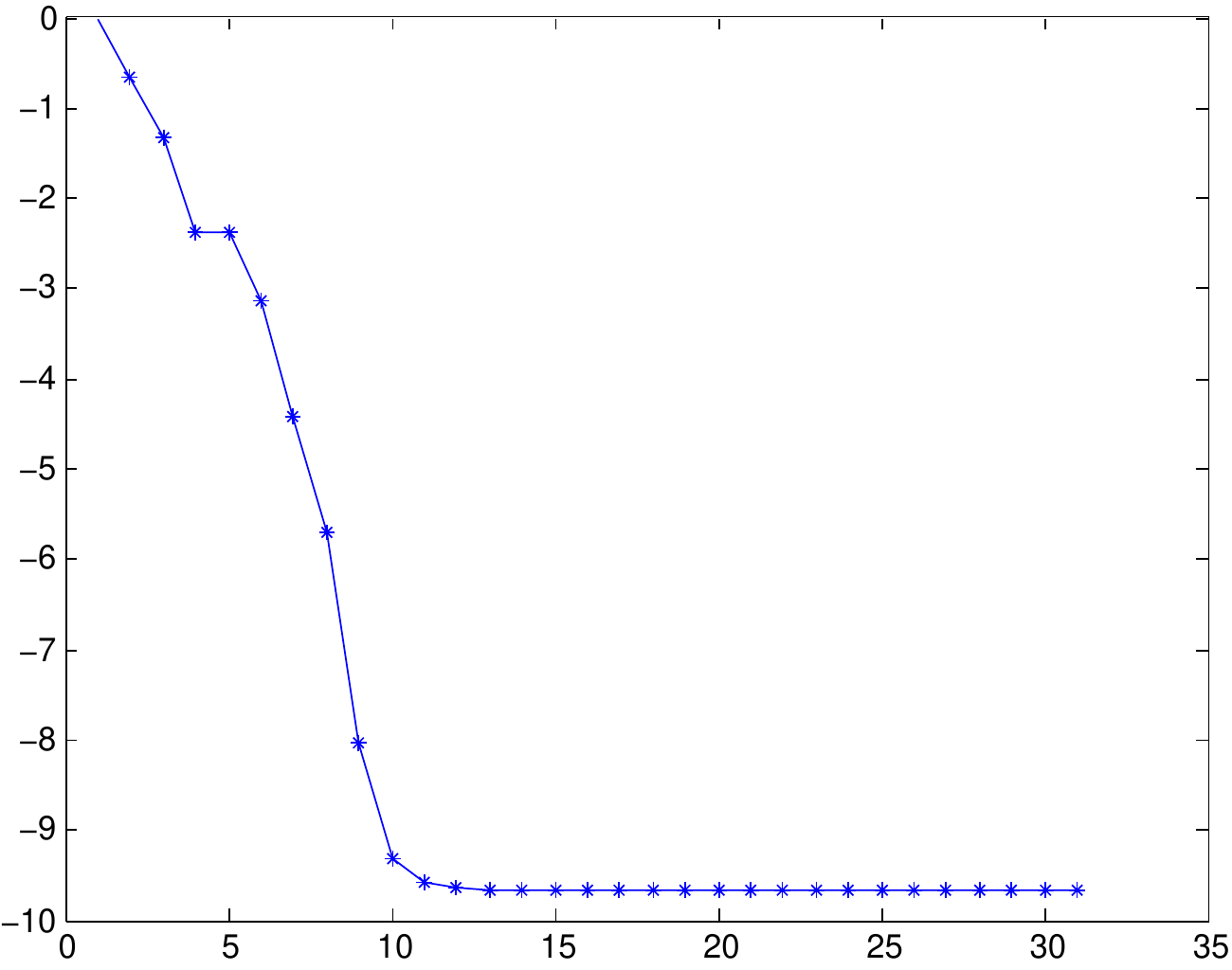}\hskip 0.6cm
\caption{Evolution of the objective functions (normalized with its starting value and in logarithmic scale) in BFGS Newton iteration for the three reconstructions in Experiment II.}
\label{FIG:Exp Rec Con His}
\end{figure}

\paragraph{Experiment II.} In the second set of numerical experiments, we consider the reconstruction of a long absorbing object. The scattering coefficient is again set to be $8$. The singular values of the matrix $\bA$ is shown in second plot of Fig.~\ref{FIG:SVD}. Note that even though the domain and the scattering coefficient in this case are the same as those in Experiment I, the finite volume mesh are different in the two cases. That is why we observe a slight difference between the singular values of the two cases. The reconstruction results are presented in Fig.~\ref{FIG:Exp Rec}. Shown, from left to right, are the true $\sigma_a$ and reconstructions with data of type (i), (ii) and (iii) respectively. The relative $L^2$ error for reconstructions are $3.24\%$, $5.68\%$, and $10.46\%$ respectively. These are comparable to the reconstructions using more expensive reconstruction methods, such as those in~\cite{ReBaHi-SIAM06}. Convergence history of the reconstructions are shown in Fig.~\ref{FIG:Exp Rec Con His}. We observed again that when the noise level is high, the algorithm converges very slow in later iterations. This is a clear indication that in these cases, we can choose a smaller $L$ value.

\begin{figure}[ht]
\centering
\includegraphics[angle=0,width=0.230\textwidth]{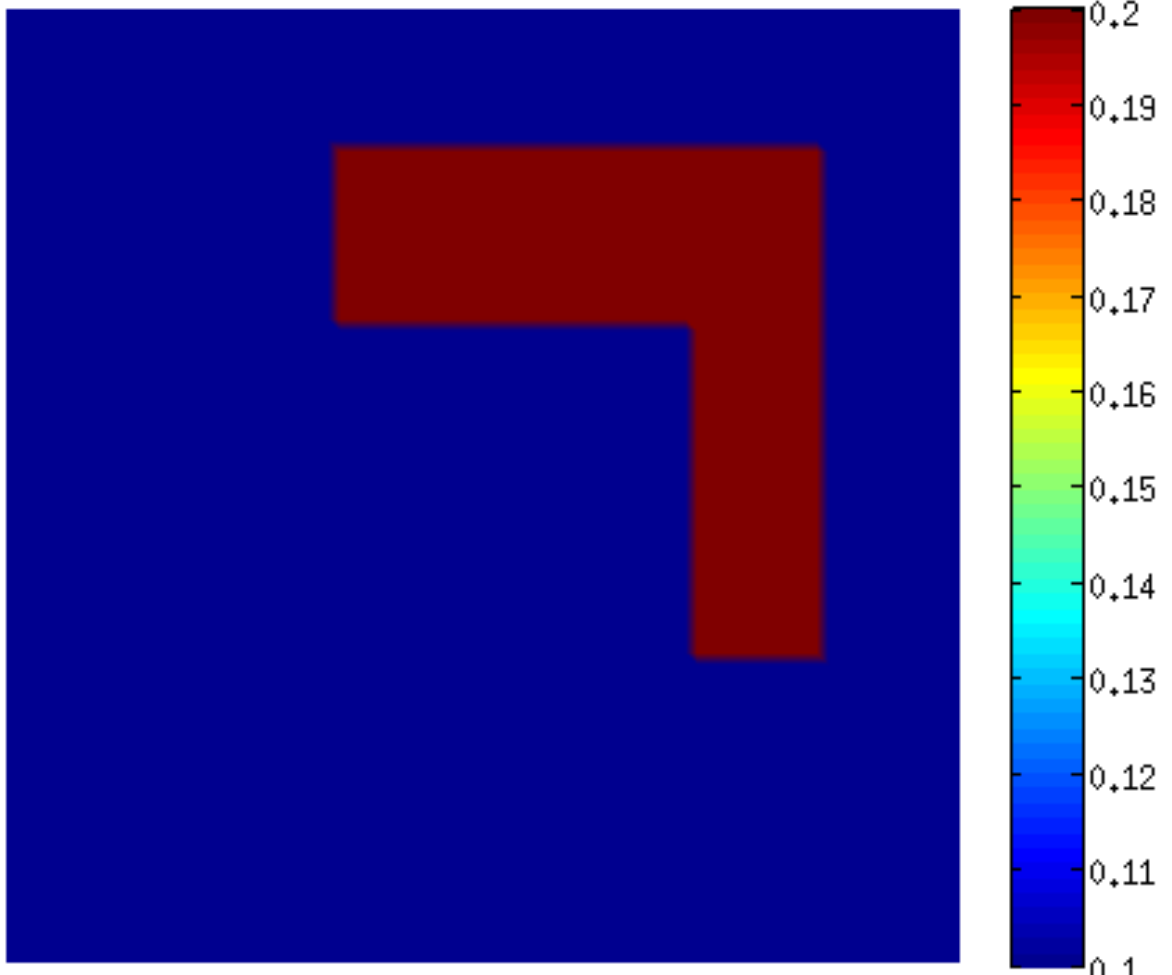}
\includegraphics[angle=0,width=0.230\textwidth]{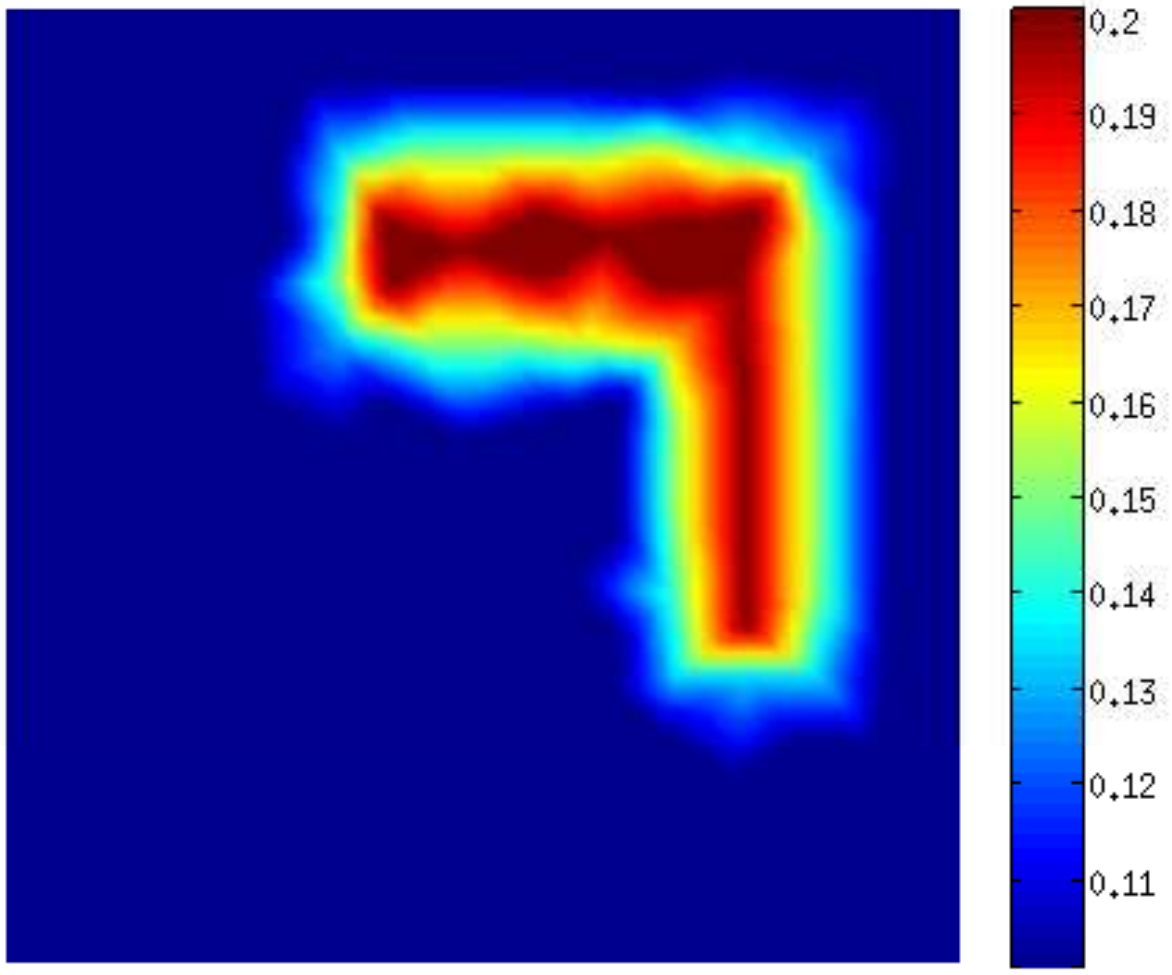}
\includegraphics[angle=0,width=0.230\textwidth]{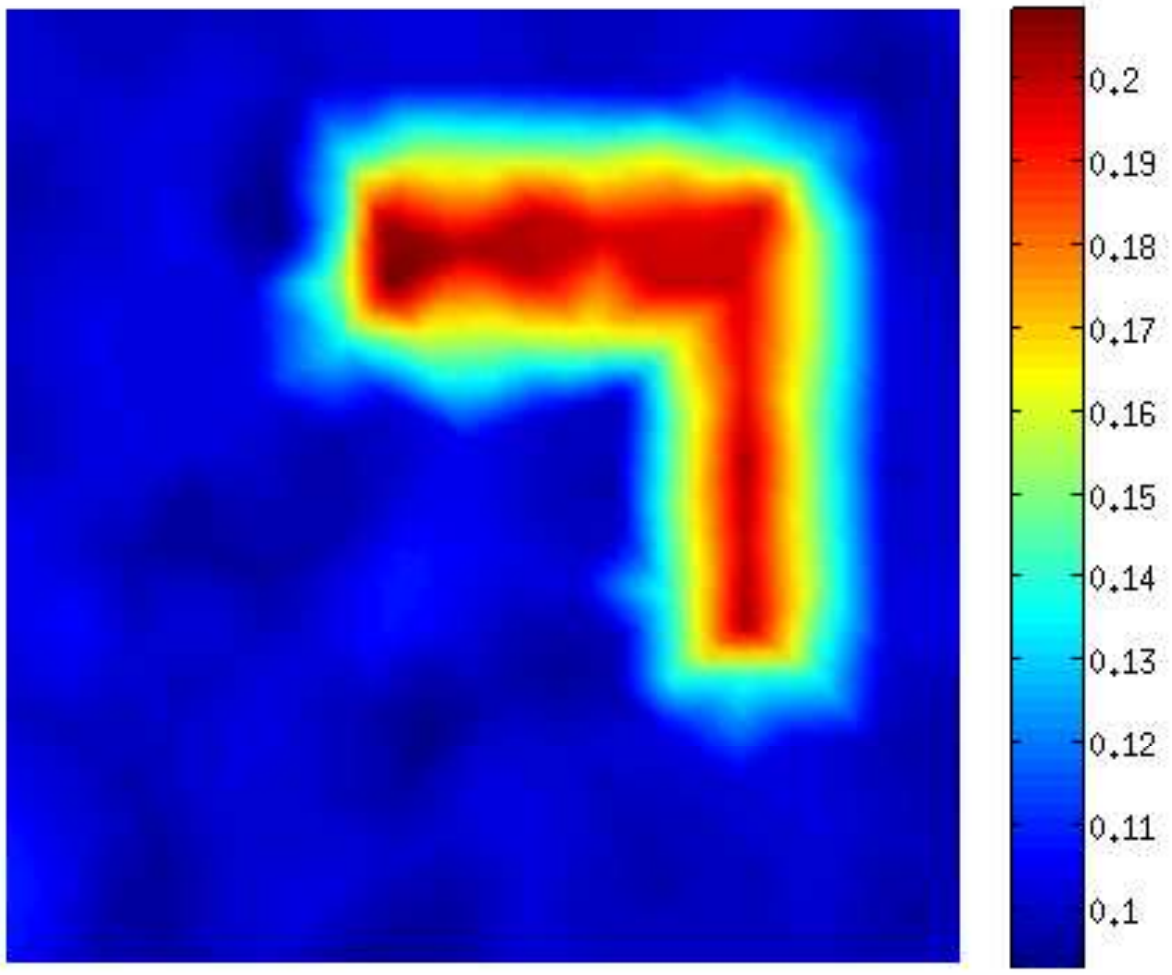}
\includegraphics[angle=0,width=0.230\textwidth]{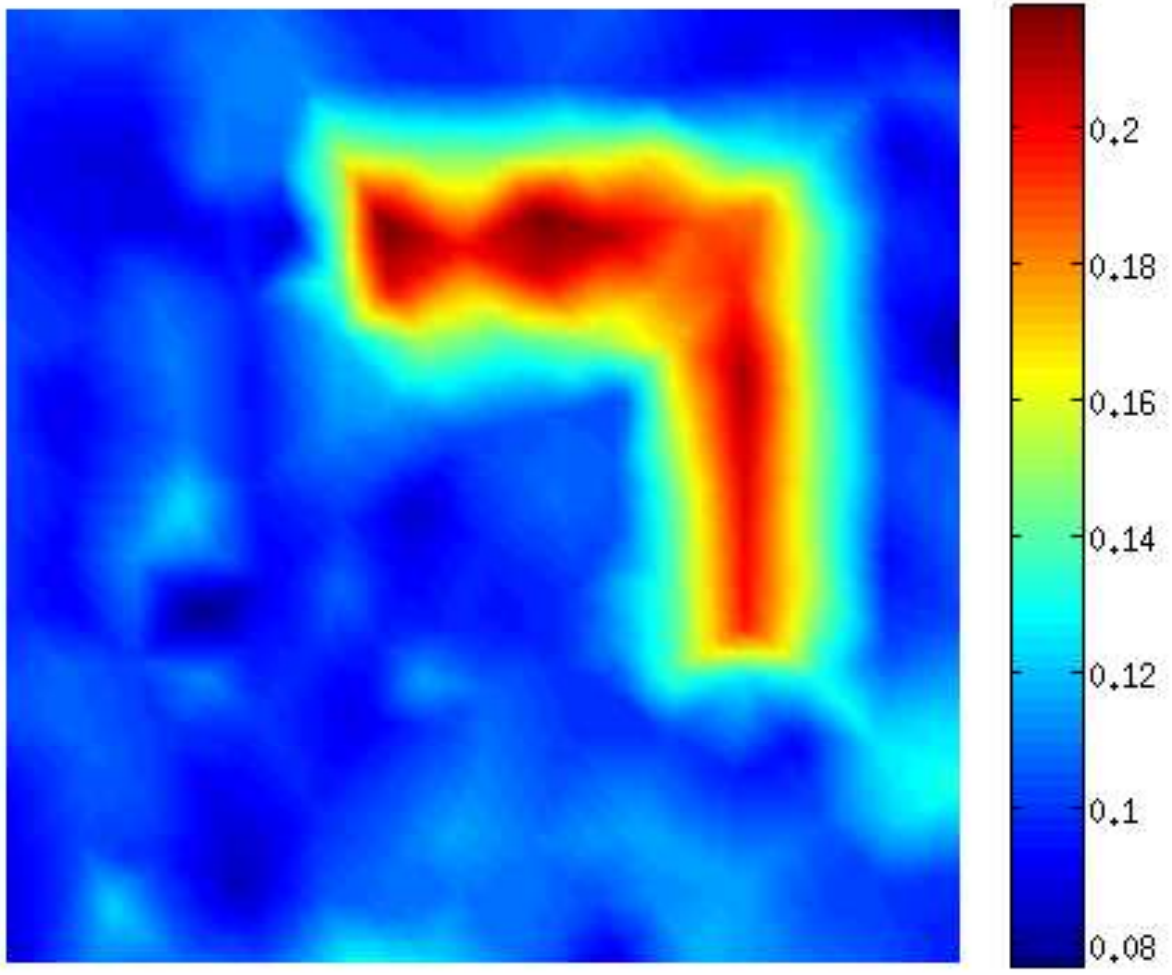}
\caption{From left to right: True absorption coefficient and absorption coefficient reconstructed with data of type (i), type (ii) and type (iii) respectively for Experiment III.}
\label{FIG:Exp Lshaped}
\end{figure}
\begin{figure}[ht]
\centering
\includegraphics[angle=0,width=0.30\textwidth]{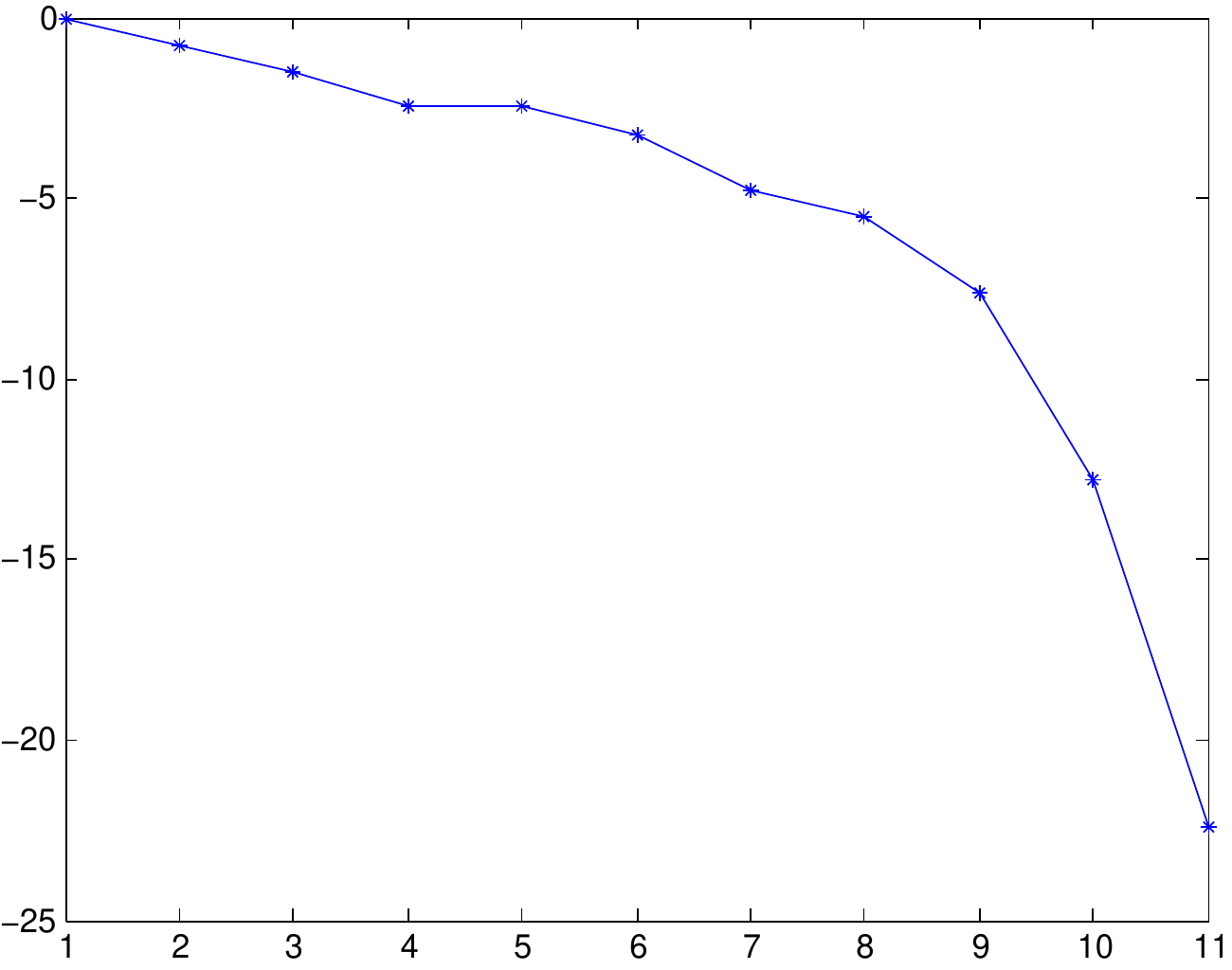}\hskip 0.6cm
\includegraphics[angle=0,width=0.30\textwidth]{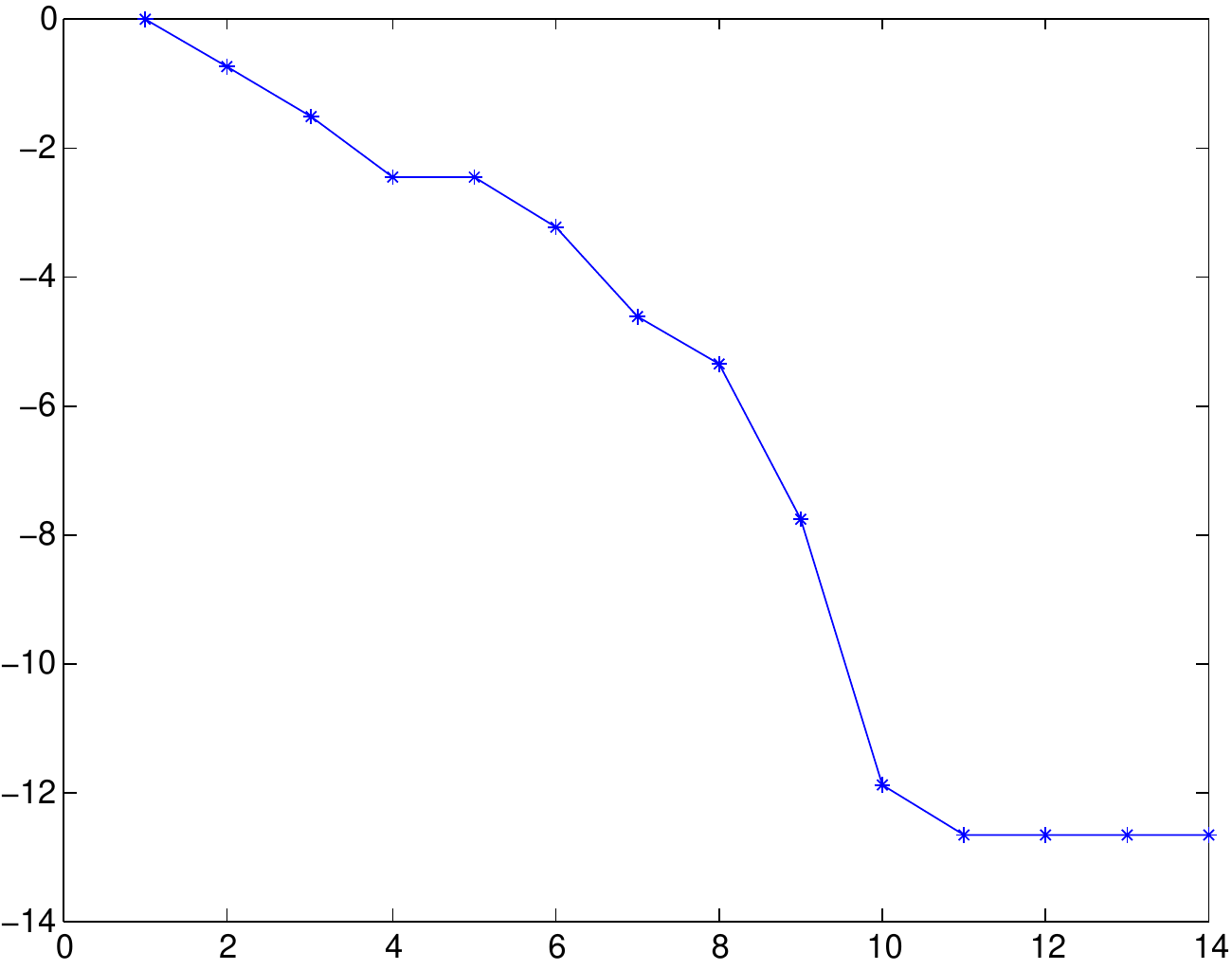}\hskip 0.6cm
\includegraphics[angle=0,width=0.30\textwidth]{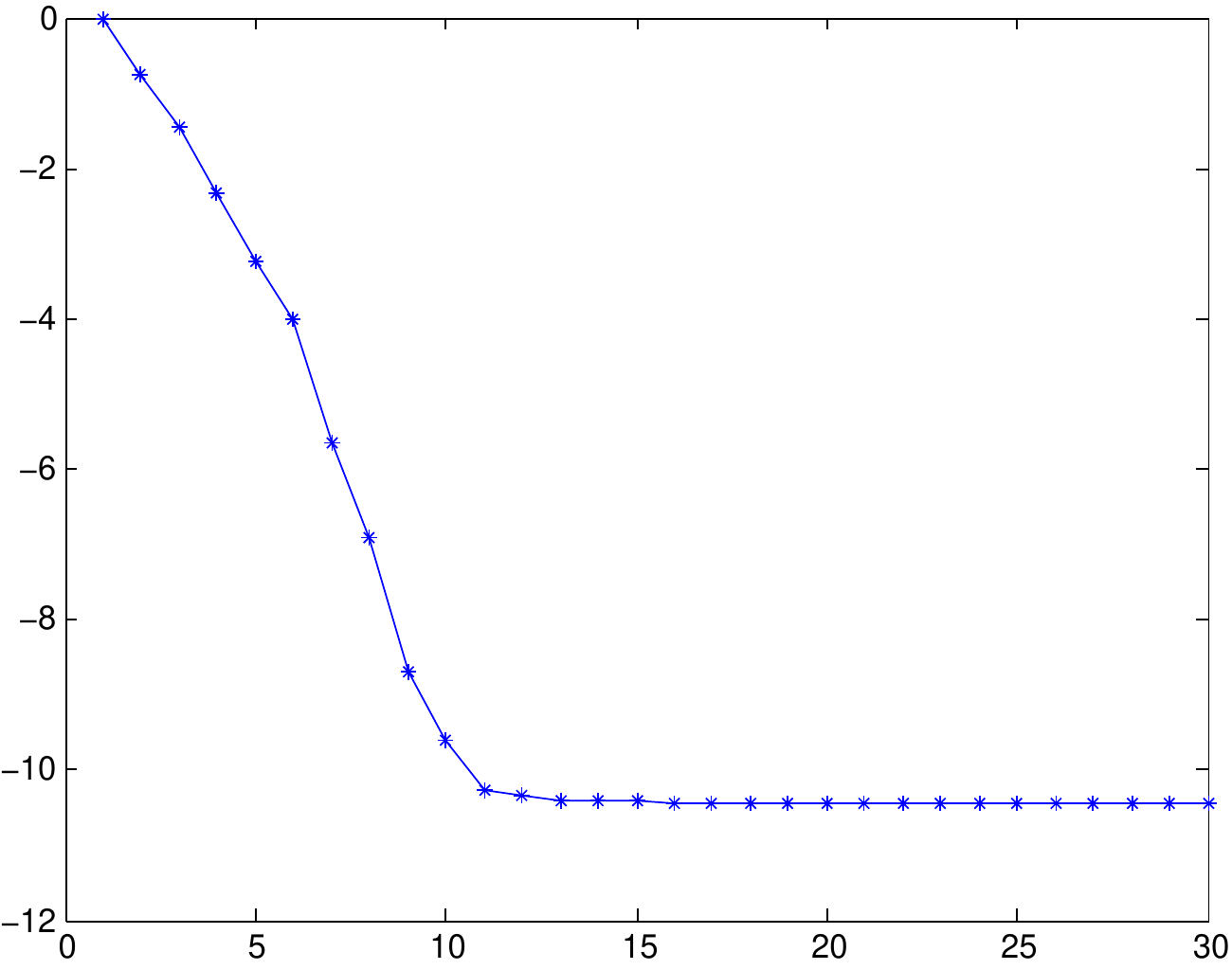}\hskip 0.6cm
\caption{Evolution of the objective functions (normalized with its starting value and in logarithmic scale) in BFGS Newton iteration for the three reconstructions in Experiment III.}
\label{FIG:Exp Lshaped Con His}
\end{figure}

\paragraph{Experiment III.} We repeat the reconstruction process for a more complicated absorption coefficient shown in the left plot of Fig.~\ref{FIG:Exp Lshaped}. The reconstruction results with data of type (i), (ii) and (iii) respectively are presented in the right plots of Fig.~\ref{FIG:Exp Lshaped}. The quality of the reconstructions are similar to those in the first two cases. The relative $L^2$ error for reconstructions in the top row are $3.20\%$, $6.68\%$, and $10.67\%$ respectively. Convergence history of the reconstructions are shown in Fig.~\ref{FIG:Exp Lshaped Con His}.

\begin{figure}[ht]
\centering
\includegraphics[angle=0,width=0.230\textwidth]{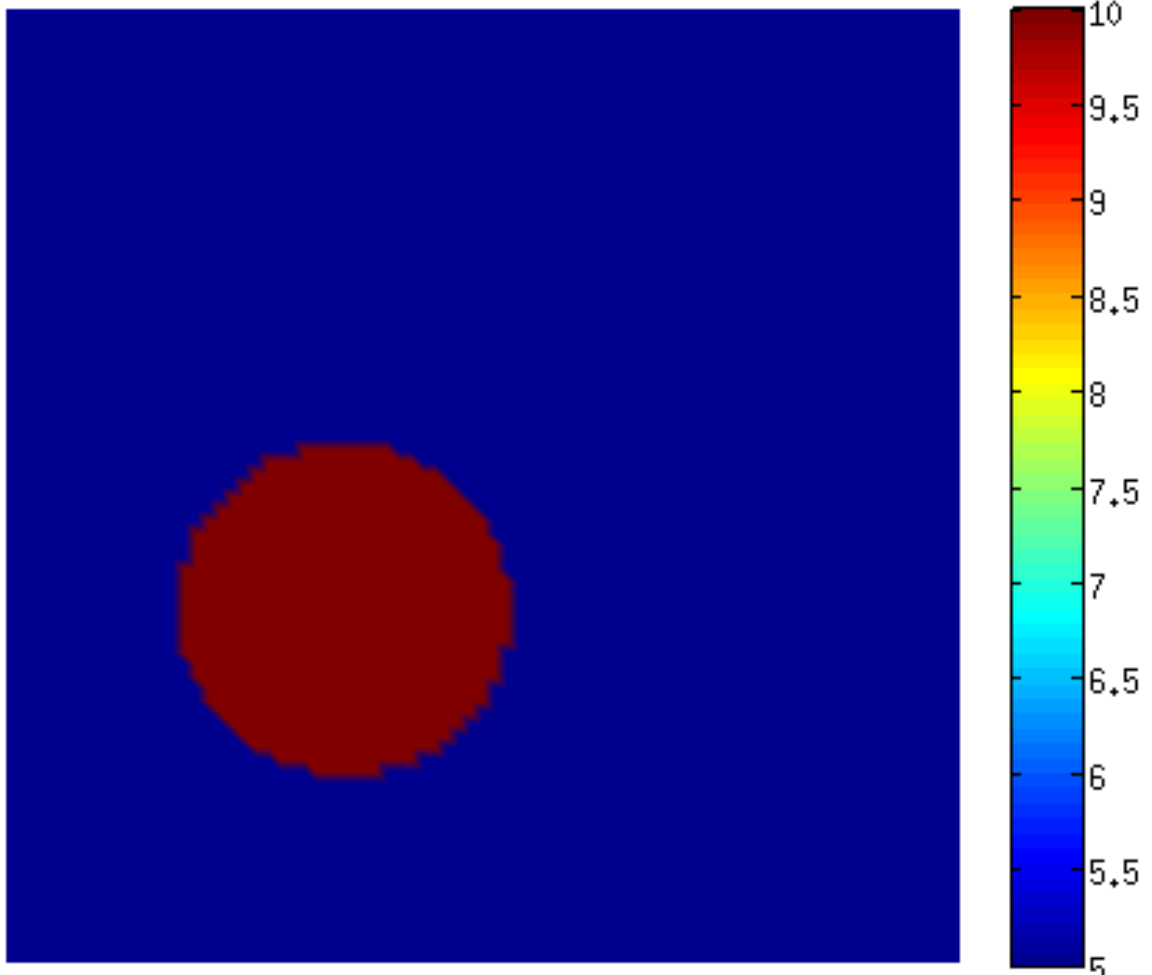}
\includegraphics[angle=0,width=0.230\textwidth]{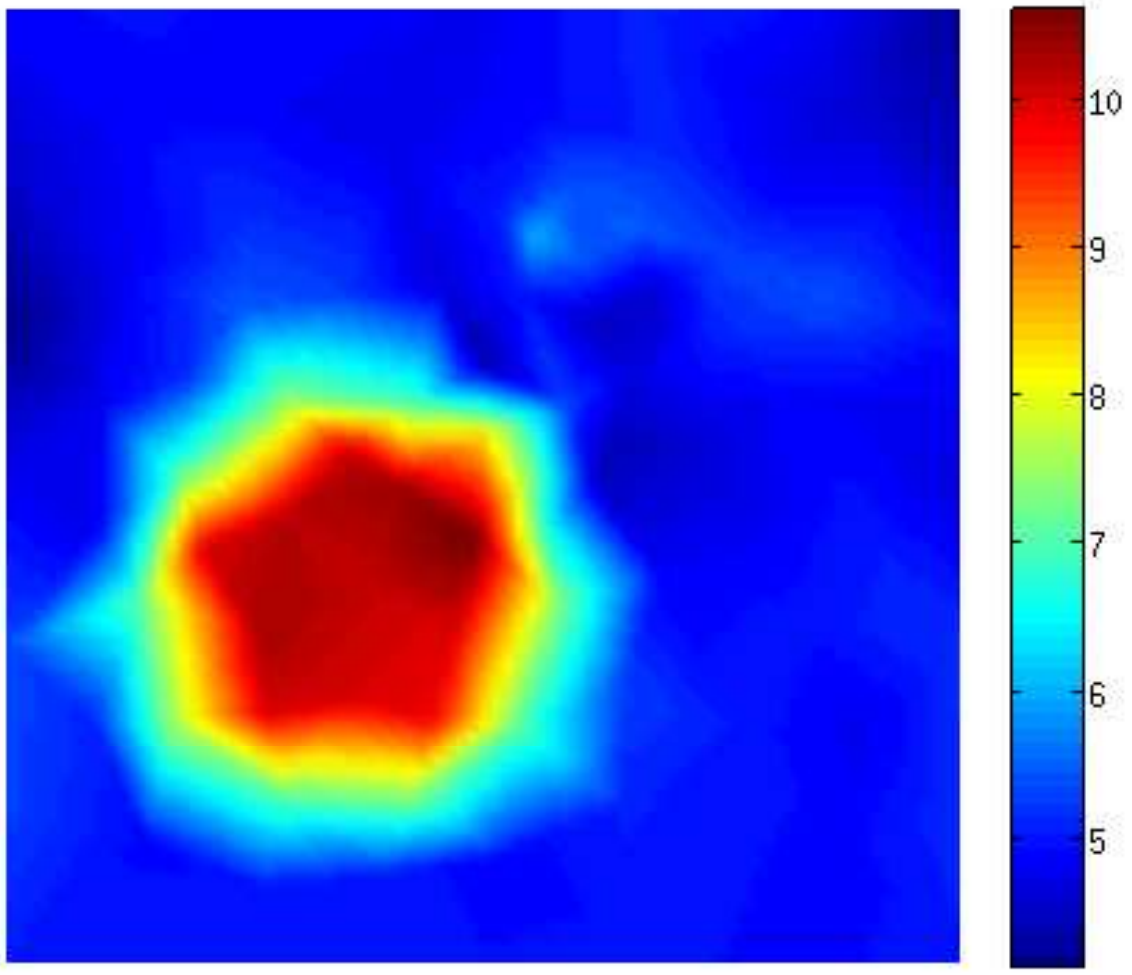}
\includegraphics[angle=0,width=0.230\textwidth]{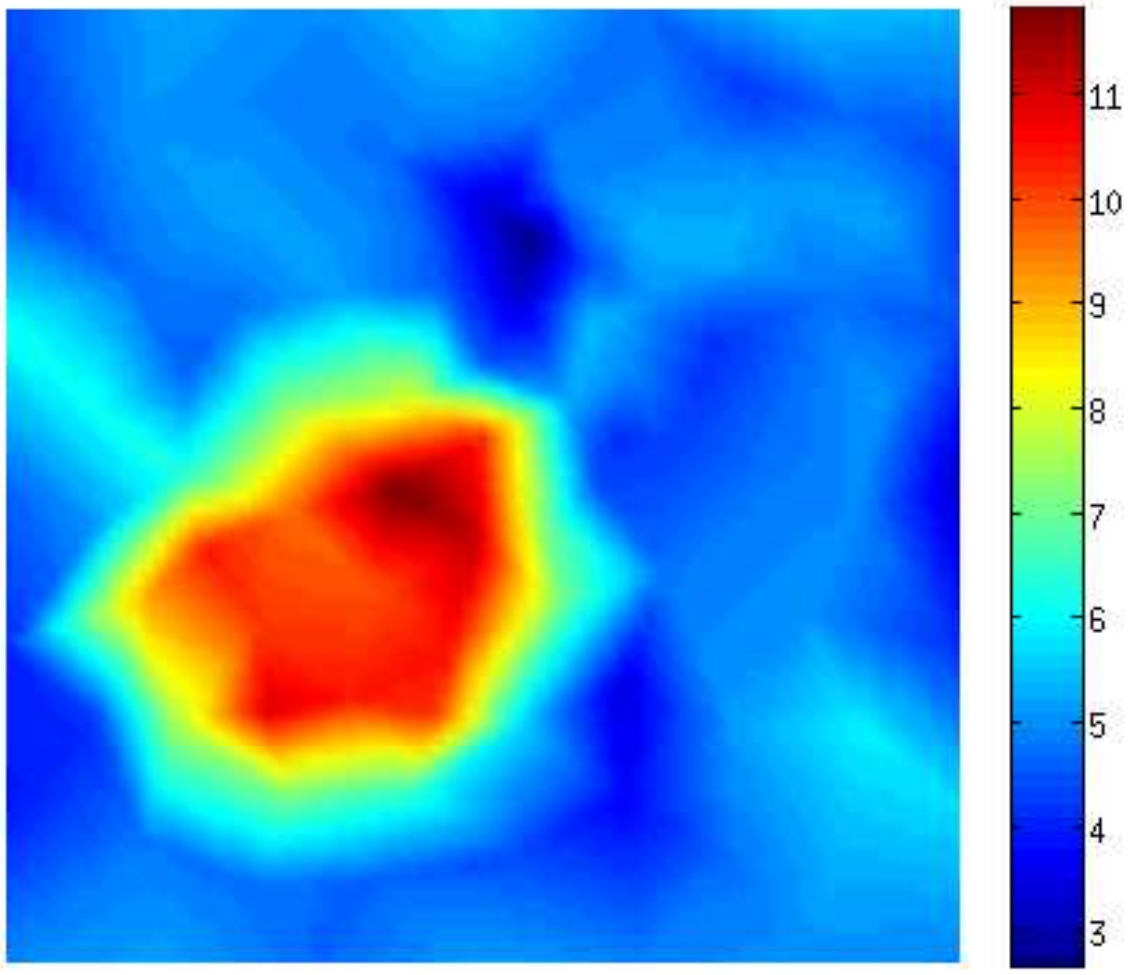}
\includegraphics[angle=0,width=0.230\textwidth]{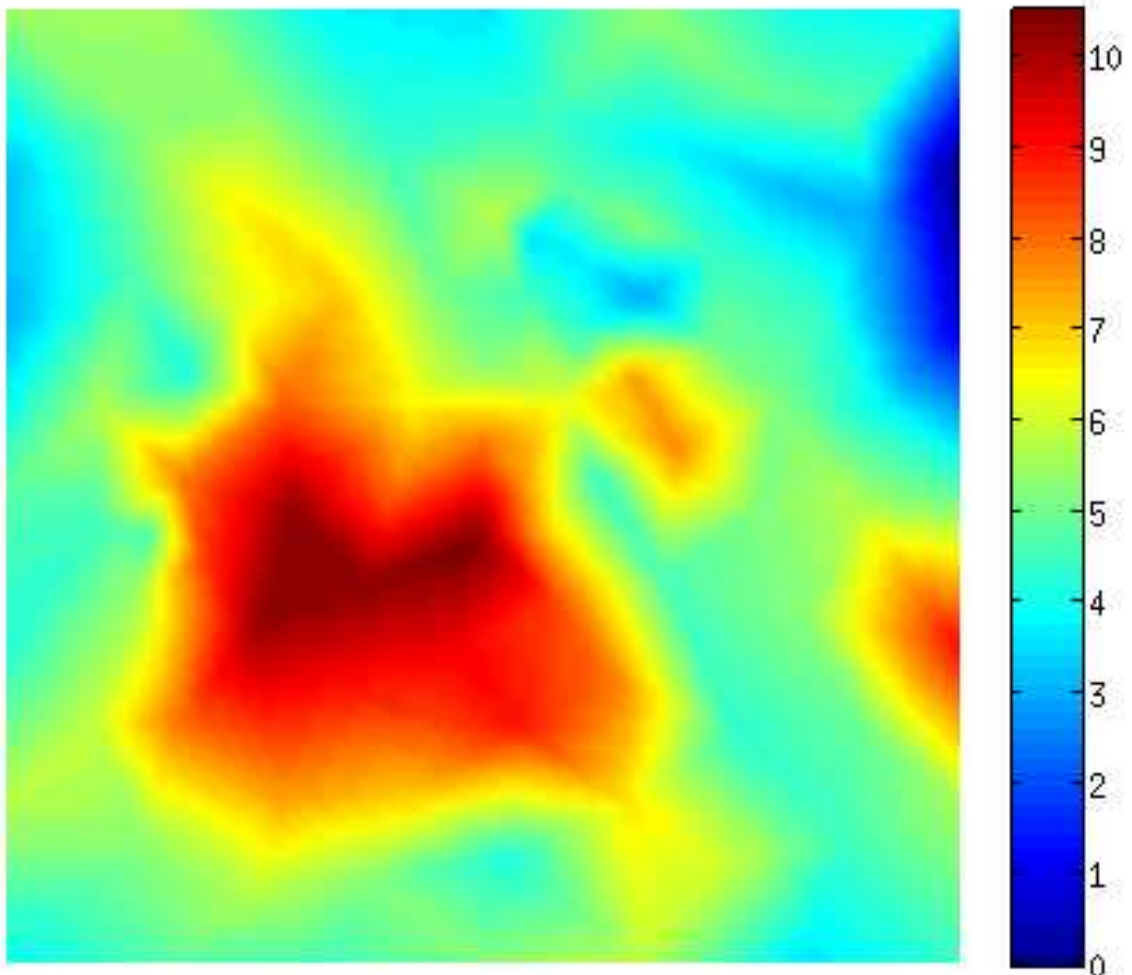}
\caption{From left to right: True scattering coefficient and scattering coefficient reconstructed with data of type (i), type (ii) and type (iii) respectively for Experiment IV.}
\label{FIG:Exp 2Cir}
\end{figure}
\paragraph{Experiment IV.} In the last numerical experiment we reconstruct the scattering coefficient assuming that the absorption coefficient is known. The absorption coefficient is assumed to be $0.2$ in a disk of radius $0.3$ centered at $(1.3, 1.4)$ and $0.1$ everywhere else. The decay rate of for large singular values is much larger than from those in Experiments I$\sim$ III as shown in the right plot of Fig.~\ref{FIG:SVD}. The true coefficient and the reconstructions with three different data are shown in Fig.~\ref{FIG:Exp 2Cir}. The relative $L^2$ error in the three reconstructions are $6.34\%$, $8.82$, $13.77\%$ respectively. The reconstructions are slightly worse than those on the absorption coefficients in the previous experiments. This is not due to the algorithm itself but due to the fact that the scattering coefficient is harder to reconstruct than the absorption coefficient, as is well-known in the diffuse optical tomography community and reflected partially in the fast decay rate of the singular values of $\bA$.

The computational costs of the reconstructions in this section are negligible after the singular value decompositions have been constructed, on the order of minutes on a reasonable desktop, such as a DELL OPTIPLEX 780 with Intel Core 2 Quad Q9400 with 8 GB of memory. The construction of the SVD, however, is very expensive. For the simulations we have here, each SVD cost on the order of $10$ computational hours on the same desktop. In cases when multiple reconstructions have to be done in the same configuration, which we believe is what has to be done in practical applications, our method will out-perform traditional methods since we only construct one SVD and use it for the rest of reconstructions.

\section{Concluding remarks}
\label{SEC:Concl}

To summarize, we proposed in this work a subspace-based minimization reconstruction strategy for solving inverse coefficient problems for the radiative transport equation, for applications in optical imaging techniques such as diffuse optical tomography. In this strategy, we factorize the map from the unknown coefficient to the observed data into the composition of a map from the coefficient to an intermediate variable (which is also unknown) with a map from the intermediate variable to the data. We then perform a spectral decomposition (SVD) of the second map which enable us to decompose the intermediate variable into a low-frequency part, which can be stably reconstructed, and a high-frequency part, which can not be stably reconstructed due to noise in the data. In the reconstruction, we reconstruct the low-frequency component of the intermediate variable analytically and the high-frequency component with an inexpensive minimization algorithm. We then reconstruct the coefficient by inverting the map from the coefficient to the intermediate variable using another inexpensive minimization algorithm. Numerical simulation results based on synthetic data demonstrated that this reconstruction strategy can be efficient and robust once the computationally expensive spectral decomposition have been computed off-line.

Even though we use the terminology ``subspace-based minimization'' following the work of Chen and collaborators~\cite{Chen-JOSA09,PaChZhYe-JOSA10,ZhCh-IP09}, there are several critical differences between our strategy and these in~\cite{Chen-JOSA09,PaChZhYe-JOSA10,ZhCh-IP09} as we have emphasized in the presentation of Section~\ref{SEC:SOM}. The main difference is the reduction of unknowns in the reconstruction following the philosophy in one of our previous work~\cite{GuReMaHi-JBO10}. In fact, it would be very interesting to combine the current algorithm, in which the intermediate variable is parameterized, with the algorithm in~\cite{GuReMaHi-JBO10}, in which the unknown coefficient is parameterized. There are also strong connections, in terms of algorithm philosophy, between our algorithm, the algorithms based on optimal grid and networks~\cite{BoDr-IP02,BoDrGu-IP08,BoDrMa-IP10} in which the measurement setup of the problem is used to determine an optimal parameterization of the unknowns for the reconstruction algorithm, and the algorithms based on sparsity~\cite{GaCaShZh-PMB11,GaYuOsWa-IP11}. Essentially, all the aforementioned strategies share the same philosophy, that is, for these severely ill-conditioned inverse problems, only low-frequency contents in the unknowns can be stably reconstructed when no extra \emph{a priori} information are available. It is thus more efficient to simply attempt to reconstruct these low-frequency contents parameterized under a good basis. The difference between the algorithms lies in their strategies to obtain that good basis.

\section*{Acknowledgement}

We would like to thank the anonymous referees for their useful comments that help us improve the quality of this paper. We would also like to thank Professor Xudong Chen (National University of Singapore) for useful discussions on the subspace-based minimization methods. This work is partially supported by National Science Foundation (NSF) through grants DMS-0914825 and DMS-1321018.



\begin{thebibliography}{10}

\bibitem{AbReHi-IP05}
{\sc G.~S. Abdoulaev, K.~Ren, and A.~H. Hielscher}, {\em Optical tomography as
  a {PDE}-constrained optimization problem}, Inverse Problems, 21 (2005),
  pp.~1507--1530.

\bibitem{Agoshkov-Book98}
{\sc V.~Agoshkov}, {\em Boundary Value Problems for the Transport Equations},
  Birkhauser, Boston, 1998.

\bibitem{ArBaLuGr-MMMTT91}
{\sc R.~Aronson, R.~L. Barbour, J.~Lubowsky, and H.~Graber}, {\em Application
  of transport theory to infra-red medical imaging}, in Modern Mathematical
  Methods in Transport Theory, W.~Greenberg and J.~Polewczak, eds.,
  Birkh\"{a}user, 1991.

\bibitem{Arridge-IP99}
{\sc S.~R. Arridge}, {\em Optical tomography in medical imaging}, Inverse
  Probl., 15 (1999), pp.~R41--R93.

\bibitem{Arridge-IP06}
{\sc S.~R. Arridge, O.~Dorn, J.~P. Kaipio, V.~Kolehmainen, M.~Schweiger,
  T.~Tarvainen, M.~Vauhkonen, and A.~Zacharopoulos}, {\em Reconstruction of
  subdomain boundaries of piecewise constant coefficients of the radiative
  transfer equation from optical tomography data}, Inverse Problems, 22 (2006),
  pp.~2175--2196.

\bibitem{Babovsky-SIAM91}
{\sc H.~Babovsky}, {\em Identification of scattering media from reflected
  flows}, SIAM J. Appl. Math., 51 (1991), pp.~1676--1704.

\bibitem{Babovsky-IP95}
\leavevmode\vrule height 2pt depth -1.6pt width 23pt, {\em An inverse model
  problem in kinetic theory}, Inverse Problems, 11 (1995), pp.~555--570.

\bibitem{Bal-IP09}
{\sc G.~Bal}, {\em Inverse transport theory and applications}, Inverse
  Problems, 25 (2009).
\newblock 053001.

\bibitem{BaLaMo-IPI08}
{\sc G.~Bal, I.~Langmore, and F.~Monard}, {\em Inverse transport with isotropic
  sources and angularly averaged measurements}, Inverse Problems and Imaging, 2
  (2008), pp.~23--42.

\bibitem{BaMo-JCP10}
{\sc G.~Bal and F.~Monard}, {\em An accurate solver for forward and inverse
  transport}, J. Comp. Phys., 229 (2010), pp.~4952--4979.

\bibitem{BaRe-IP05}
{\sc G.~Bal and K.~Ren}, {\em Atmospheric concentration profile reconstructions
  from radiation measurements}, Inverse Problems, 21 (2005), pp.~153--168.

\bibitem{BaTa-SIAM07}
{\sc G.~Bal and A.~Tamasan}, {\em Inverse source problems in transport
  equations}, SIAM J. Math. Anal., 39 (2007), pp.~57--76.

\bibitem{Barichello-Pro02}
{\sc L.~B. Barichello}, {\em Some comments concerning inverse problems in
  particle transport theory}, in Proceedings of 4th International Conference on
  Inverse Problems in Engineering, Rio de Janeiro, Brazil, 2002.

\bibitem{Belleni-Morante-TTSP03}
{\sc A.~Belleni-Morante}, {\em An inverse problem for photon transport in
  interstellar clouds}, Transport Theory and Statistical Physics, 32 (2003),
  pp.~73--91.

\bibitem{BeRo-JMAA02}
{\sc A.~Belleni-Morante and G.~F. Roach}, {\em Gamma ray transport in the
  cardiac region: an inverse problem}, J. Math. Anal. Appl., 269 (2002),
  pp.~200--215.

\bibitem{Boman-Thesis07}
{\sc E.~Boman}, {\em Radiotherapy Forward and Inverse Problem Applying
  Boltzmann Transport Equations}, PhD thesis, University of Kuopio, Filand,
  Kuopio, Filand, 2007.

\bibitem{BoDr-IP02}
{\sc L.~Borcea and V.~Druskin}, {\em Optimal finite difference grids for direct
  and inverse {Sturm-Liouville} problems}, Inverse Probl., 18 (2002),
  pp.~979--1001.

\bibitem{BoDrGu-IP08}
{\sc L.~Borcea, V.~Druskin, and F.~Guevara~Vasquez}, {\em Electrical impedance
  tomography with resistor networks}, Inverse Problems, 24 (2008).
\newblock 035013.

\bibitem{BoDrMa-IP10}
{\sc L.~Borcea, V.~Druskin, and A.~V. Mamonov}, {\em Circular resistor networks
  for electrical impedance tomography with partial boundary measurements},
  Inverse Problems, 26 (2010).
\newblock 045010.

\bibitem{CaXuAl-IEEE03}
{\sc W.~Cai, M.~Xu, and R.~R. Alfano}, {\em Three-dimensional radiative
  transfer tomography for turbid media}, IEEE Journal of Selected Topics in
  Quantum Electronics, 9 (2003), pp.~189--198.

\bibitem{Case-PF73}
{\sc K.~M. Case}, {\em Inverse problem in transport theory}, Phys. Fluids, 16
  (1973), pp.~1607--1611.

\bibitem{Case-PF75}
\leavevmode\vrule height 2pt depth -1.6pt width 23pt, {\em Inverse problem in
  transport theory. {II}}, Phys. Fluids, 18 (1975), pp.~927--930.

\bibitem{Case-PF77}
\leavevmode\vrule height 2pt depth -1.6pt width 23pt, {\em Inverse problem in
  transport theory. {III}}, Phys. Fluids, 20 (1977), pp.~2031--2036.

\bibitem{Chahine-JAS70}
{\sc M.~T. Chahine}, {\em Inverse problems in radiative transfer: Determination
  of atmospheric parameters}, J. Atmos. Sci., 27 (1970), pp.~960--967.

\bibitem{Chen-JOSA09}
{\sc X.~Chen}, {\em Application of signal-subspace and optimization methods in
  reconstructing extended scatterers}, J. Opt. Soc. Am. A, 26 (2009),
  pp.~1022--1026.

\bibitem{ChSt-CPDE96}
{\sc M.~Choulli and P.~Stefanov}, {\em Inverse scattering and inverse boundary
  value problem for the linear {Boltzmann} equation}, Comm. Part. Diff. Eqn.,
  21 (1996), pp.~763--785.

\bibitem{ChSt-IP96}
\leavevmode\vrule height 2pt depth -1.6pt width 23pt, {\em Reconstruction of
  the coefficients of the stationary transport equation from boundary
  measurements}, Inverse Probl., 12 (1996), pp.~L19--L23.

\bibitem{DeVo-JCP02}
{\sc A.~Dedner and P.~Vollm\"oller}, {\em An adaptive higher order method for
  solving the radiation transport equation on unstructured grids}, J. Comput.
  Phys., 178 (2002), pp.~263--289.

\bibitem{Dorn-CAMQ01}
{\sc O.~Dorn}, {\em Shape reconstruction in scattering media with void using a
  transport model and level sets}, Canad. Appl. Math. Quart., 10 (2001),
  pp.~239--275.

\bibitem{DrBa-M2AS05}
{\sc F.~Dragoni and L.~Barletti}, {\em An inverse problem for two-frequency
  photon transport in a slab}, Math. Meth. Appl. Sci., 28 (2005),
  pp.~1695--1714.

\bibitem{DuKl-JCP02}
{\sc B.~Dubroca and A.~Klar}, {\em Half-moment closure for radiative transfer
  equations}, J. Comput. Phys., 180 (2002), pp.~584--596.

\bibitem{ElDuMcEm-JOSA88}
{\sc R.~A. Elliott, T.~Duracz, N.~J. McCormick, and D.~R. Emmons}, {\em
  Experimental test of a time-dependent inverse radiative transfer algorithm
  for estimating scattering parameters}, J. Opt. Soc. Am. A, 5 (1988),
  pp.~366--373.

\bibitem{GaCaShZh-PMB11}
{\sc H.~Gao, J.~F. Cai, Z.~Shen, and H.~Zhao}, {\em Robust principal component
  analysis-based four-dimensional computed tomography}, Phys. Med. Bio., 56
  (2011), pp.~3181--3198.

\bibitem{GaYuOsWa-IP11}
{\sc H.~Gao, H.~Yu, S.~Osher, and G.~Wang}, {\em Multi-energy {CT} based on a
  prior rank, intensity and sparsity model ({PRISM})}, Inverse Problems, 27
  (2011).
\newblock 115012.

\bibitem{GaZh-TTSP09}
{\sc H.~Gao and H.~Zhao}, {\em A fast forward solver of radiative transfer
  equation}, Transport Theory and Statistical Physics, 38 (2009), pp.~149--192.

\bibitem{GaZh-OE10A}
\leavevmode\vrule height 2pt depth -1.6pt width 23pt, {\em Multilevel
  bioluminescence tomography based on radiative transfer equation. part 1: l1
  regularization}, Optics Express, 18 (2010), pp.~1854--1871.

\bibitem{GaZh-OE10B}
\leavevmode\vrule height 2pt depth -1.6pt width 23pt, {\em Multilevel
  bioluminescence tomography based on radiative transfer equation. part 2:
  total variation and l1 data fidelity}, Optics Express, 18 (2010),
  pp.~2894--2912.

\bibitem{GiHoWaFaMo-PMB06}
{\sc K.~A. Gifford, J.~L. Horton, T.~A. Wareing, G.~Failla, and F.~Mourtada},
  {\em Comparison of a finite-element multigroup discrete-ordinates code with
  monte carlo for radiotherapy calculations}, Phys. Med. Biol., 51 (2006),
  p.~2253–2265.

\bibitem{GoKi-IP09}
{\sc P.~Gonz\'alez-Rodr\'iguez and A.~D. Kim}, {\em Reflectance optical
  tomography in epithelial tissues}, Inverse Problems, 25 (2009).
\newblock 015001.

\bibitem{GuReHi-AO07}
{\sc X.~Gu, K.~Ren, and A.~H. Hielscher}, {\em Frequency-domain sensitivity
  analysis for small imaging domains using the equation of radiative transfer},
  Applied Optics, 46 (2007), pp.~6669--6679.

\bibitem{GuReMaHi-JBO10}
{\sc X.~Gu, K.~Ren, J.~Masciotti, and A.~H. Hielscher}, {\em Parametric image
  reconstruction using the discrete cosine transform for optical tomography},
  J. Biomed. Opt., 14 (2010).
\newblock Art. No. 064003.

\bibitem{HeGr-AJ41}
{\sc L.~G. Henyey and J.~L. Greenstein}, {\em Diffuse radiation in the galaxy},
  Astrophys. J., 90 (1941), pp.~70--83.

\bibitem{JiScPaJi-PMB12}
{\sc X.~Jia, J.~Sch\"umann, H.~Paganetti, and S.~B. Jiang}, {\em Gpu-based fast
  monte carlo dose calculation for proton therapy}, Phys. Med. Biol., 57
  (2012), pp.~7783--7797.

\bibitem{KaDa-TTSP79}
{\sc M.~Kanal and J.~A. Davies}, {\em A multidimensional inverse problem in
  transport theory}, Transport Theory Statist. Phys., 8 (1979), pp.~99--115.

\bibitem{KaMo-JMP78B}
{\sc M.~Kanal and H.~E. Moses}, {\em Direct-inverse problems in transport
  theory. {I}. {The} inverse problem}, J. Math. Phys., 19 (1978),
  pp.~1793--1798.

\bibitem{KaMo-JMP78}
\leavevmode\vrule height 2pt depth -1.6pt width 23pt, {\em Direct-inverse
  problems in transport theory, the inverse albedo problem for a finite
  medium}, J. Math. Phys., 19 (1978), pp.~2641--2645.

\bibitem{KiIs-JCP99}
{\sc A.~D. Kim and A.~Ishimaru}, {\em Chebyshev spectral method for radiative
  transfer equations applied to electromagnetic wave propagation and scattering
  in a discrete random medium}, J. Comput. Phys., 152 (1999), pp.~264--280.

\bibitem{KiMo-JCP03}
{\sc A.~D. Kim and M.~Moscoso}, {\em Radiative transfer computations for
  optical beams}, J. Comput. Phys., 185 (2003), pp.~50--60.

\bibitem{KiMo-IP06}
\leavevmode\vrule height 2pt depth -1.6pt width 23pt, {\em Radiative transport
  theory for optical molecular imaging}, Inverse Problems, 22 (2006),
  pp.~23--42.

\bibitem{KiFlYaKaHi-BOE10}
{\sc H.~K. Kim, M.~Flexman, D.~Y. Yamashiro, J.~Kandel, and A.~H. Hielscher},
  {\em {PDE}-constrained multispectral imaging of tissue chromophores with the
  equation of radiative transfer}, Biomedical Optics Express, 1 (2010),
  pp.~812--824.

\bibitem{KiHi-IP09}
{\sc H.~K. Kim and A.~H. Hielscher}, {\em A {PDE-constrained} {SQP} algorithm
  for optical tomography based on the frequency-domain equation of radiative
  transfer}, Inverse Problems, 25 (2009).
\newblock 015010.

\bibitem{KlYa-SIAM07}
{\sc M.~V. Klibanov and M.~Yamamoto}, {\em Exact controllability of the time
  dependent transport equation}, SIAM J. Control Optim., 46 (2007),
  pp.~2071--2095.

\bibitem{KlLa-JCP06}
{\sc A.~D. Klose and E.~W. Larsen}, {\em Light transport in biological tissue
  based on the simplified spherical harmonics equations}, J. Comput. Phys., 220
  (2006), pp.~441--470.

\bibitem{KlNtHi-JCP05}
{\sc A.~D. Klose, V.~Ntziachristos, and A.~H. Hielscher}, {\em The inverse
  source problem based on the radiative transfer equation in optical molecular
  imaging}, J. Comput. Phys., 202 (2005), pp.~323--345.

\bibitem{Langmore-IP08}
{\sc I.~Langmore}, {\em The stationary transport problem with angularly
  averaged measurements}, Inverse Problems, 24 (2008).
\newblock 015024.

\bibitem{Larsen-TTSP88}
{\sc E.~W. Larsen}, {\em Solution of three-dimensional inverse transport
  problems}, Transport Theory and Statistical Physics, 17 (1988), pp.~147--167.

\bibitem{LaFoGiDwHi-JBO07}
{\sc J.~M. Lasker, C.~J. Fong, D.~T. Ginat, E.~Dwyer, and A.~H. Hielscher},
  {\em Dynamic optical imaging of vascular and metabolic reactivity in
  rheumatoid joints}, J. Biomed. Optics, 12 (2007).
\newblock 052001.

\bibitem{LeMi-Book93}
{\sc E.~E. Lewis and W.~F. Miller}, {\em {Computational} {Methods} of {Neutron}
  {Transport}}, American Nuclear Society, La Grange Park, IL, 1993.

\bibitem{Li-JQSRT01}
{\sc H.-Y. Li}, {\em A two-dimensional cylindrical inverse source problem in
  radiative transfer}, J. Quant. Spectrosc. Radiat. Transfer, 69 (2001),
  pp.~403--414.

\bibitem{Luo-RPC98}
{\sc Z.~Luo}, {\em An overview of the bipartition model for charged particle
  transport}, Radiation Physics and Chemistry, 53 (1998), pp.~305--327.

\bibitem{MaRe-CMS14}
{\sc A.~V. Mamonov and K.~Ren}, {\em Quantitative photoacoustic imaging in
  radiative transport regime}, Comm. Math. Sci., 12 (2014), pp.~201--234.

\bibitem{McCormick-TTSP84}
{\sc N.~J. McCormick}, {\em Recent developments in inverse scattering transport
  methods}, Transport Theory and Statistical Physics, 13 (1984), pp.~15--28.

\bibitem{McCormick-TTSP86}
\leavevmode\vrule height 2pt depth -1.6pt width 23pt, {\em Methods for solving
  inverse problems for radiation transport - {An} update}, Transport Theory and
  Statistical Physics, 15 (1986), pp.~159--172.

\bibitem{McCormick-NSE92}
\leavevmode\vrule height 2pt depth -1.6pt width 23pt, {\em Inverse radiative
  transfer problems: {A} review}, Nuclear Science and Engineering, 112 (1992),
  pp.~185--198.

\bibitem{McCormick-JOSA04}
\leavevmode\vrule height 2pt depth -1.6pt width 23pt, {\em Analytic inverse
  radiative transfer equations for atmospheric and hydrologic optics}, J. Opt.
  Soc. Am. A, 21 (2004), pp.~1009--1017.

\bibitem{PaChZhYe-JOSA10}
{\sc L.~Pan, X.~Chen, Y.~Zhong, and S.~P. Yeo}, {\em Comparison among the
  variants of subspace-based optimization method for addressing inverse
  scattering problems: transverse electric case}, J. Opt. Soc. Am. A, 27
  (2010), pp.~2208--2215.

\bibitem{Panchenko-IP93}
{\sc A.~N. Panchenko}, {\em Inverse source problem of radiative transfer: a
  special case of the attenuated {Radon} transform}, Inverse Problems, 9
  (1993), pp.~321--337.

\bibitem{Ren-CiCP10}
{\sc K.~Ren}, {\em Recent developments in numerical techniques for
  transport-based medical imaging methods}, Commun. Comput. Phys., 8 (2010),
  pp.~1--50.

\bibitem{ReAbBaHi-OL04}
{\sc K.~Ren, G.~S. Abdoulaev, G.~Bal, and A.~H. Hielscher}, {\em Algorithm for
  solving the equation of radiative transfer in the frequency domain}, Optics
  Lett., 29 (2004), pp.~578--580.

\bibitem{ReBaHi-SIAM06}
{\sc K.~Ren, G.~Bal, and A.~H. Hielscher}, {\em Frequency domain optical
  tomography based on the equation of radiative transfer}, SIAM J. Sci.
  Comput., 28 (2006), pp.~1463--1489.

\bibitem{ReBaHi-AO07}
\leavevmode\vrule height 2pt depth -1.6pt width 23pt, {\em Transport- and
  diffusion-based optical tomography in small domains: {A} comparative study},
  Applied Optics, 46 (2007), pp.~6669--6679.

\bibitem{ScMa-IPI07}
{\sc J.~C. Schotland and V.~A. Markel}, {\em {Fourier-Laplace} structure of the
  inverse scattering problem for the radiative transport equation}, Inverse
  Problems in Imaging, 1 (2007), pp.~181--188.

\bibitem{ShFeOlMa-SIAM99}
{\sc D.~M. Shepard, M.~C. Ferris, G.~H. Olivera, and T.~R. Machie}, {\em
  Optimizing the delivery of radiation therapy to cancer patients}, SIAM Rev.,
  41 (1999), pp.~721--744.

\bibitem{SoVeStCh-TTSP04}
{\sc R.~P. Souto, H.~F.~C. Velho, S.~Stephany, and E.~S. Chalhoub}, {\em
  Performance analysis of radiative transfer algorithms for inverse hydrologic
  optics in a parallel environment}, Transport Theory and Statistical Physics,
  33 (2004), pp.~449--468.

\bibitem{StUh-APDE08}
{\sc P.~Stefanov and G.~Uhlmann}, {\em An inverse source problem in optical
  molecular imaging}, Anal. PDE, 1 (2008), pp.~115--126.

\bibitem{Tamasan-IP02}
{\sc A.~Tamasan}, {\em An inverse boundary value problem in two-dimensional
  transport}, Inverse Problems, 18 (2002), pp.~209--219.

\bibitem{Tamasan-CM03}
\leavevmode\vrule height 2pt depth -1.6pt width 23pt, {\em Optical tomography
  in weakly anisotropic scattering media}, in Contemporary Mathematics, AMS,
  Providence, RI, 2003.

\bibitem{TaHaHa-IP13}
{\sc J.~Tang, W.~Han, and B.~Han}, {\em A theoretical study for {RTE} based
  parameter identification problems}, Inverse Problems, 29 (2013).
\newblock 095002.

\bibitem{TaVaAr-JQSRT08}
{\sc T.~Tarvainen, M.~Vaukhonen, and S.~R. Arridge}, {\em {Gauss-Newton}
  reconstruction method for optical tomography using the finite element
  solution of the radiative transfer equation}, J. Quant. Spectrosc. Radiat.
  Transfer, 109 (2008), pp.~2767--2778.

\bibitem{TiHeToSiAlPyUl-PMB08}
{\sc L.~Tillikainen, H.~Helminen, T.~Torsti, S.~Siljam\"aki, J.~Alakuijala,
  J.~Pyyry, and W.~Ulmer}, {\em A 3d pencil-beam-based superposition algorithm
  for photon dose calculation in heterogeneous media}, Phys. Med. Biol., 53
  (2008), pp.~3821--3839.

\bibitem{TiRoNeCa-Pro02}
{\sc M.~J.~B. Tito, N.~C. Roberty, A.~J.~S. Neto, and J.~B. Cabrejo}, {\em
  Inverse radiative transfer problems in two-dimensional participating media},
  in Proceedings of 4th International Conference on Inverse Problems in
  Engineering, Rio de Janeiro, Brazil, 2002.

\bibitem{VaWaMcFaSaMo-PMB10}
{\sc O.~N. Vassiliev, T.~A. Wareing, J.~McGhee, G.~Failla, M.~Salehpour, and
  F.~Mourtada}, {\em Validation of a new grid-based boltzmann equation solver
  for dose calculation in radiotherapy with photon beams}, Phys. Med. Biol., 55
  (2010), pp.~581--598.

\bibitem{Wang-AMC00}
{\sc A.~P. Wang}, {\em Inverse problems in radiative transfer}, Appl. Math.
  Comput., 116 (2000), pp.~39--48.

\bibitem{WaUe-ASS89}
{\sc A.~P. Wang and S.~Ueno}, {\em An inverse problem in a three-dimensional
  radiative transfer}, Astrophys. Space Sci., 155 (1989), pp.~105--111.

\bibitem{Wang-AIHP99}
{\sc J.-N. Wang}, {\em Stability estimates of an inverse problem for the
  stationary transport equation}, Ann. Inst. Henri Poincar\'e, 70 (1999),
  pp.~473--495.

\bibitem{YiWuSu-JCP96}
{\sc J.~Ying, F.~Wu, and W.~Sun}, {\em Simultaneous reconstruction of two
  parameters for transport equation in a stratified half-space}, J. Comput.
  Phys., 125 (1996), pp.~434--439.

\bibitem{ZhCh-IP09}
{\sc Y.~Zhong and X.~Chen}, {\em Twofold subspace-based optimization method for
  solving inverse scattering problems}, Inverse Problems, 25 (2009).
\newblock 085003.

\end{thebibliography}
\end{document}